\documentclass[11pt]{article}
\usepackage{amsfonts}
\usepackage{color}
\newtheorem{theo}{Theorem}[section]
\newtheorem{lem}[theo]{Lemma}

\newtheorem{defi}{Definition}[section]

\newcommand{\mysection}[1]{\section{#1} \setcounter{equation}{0}}
\newcommand{\proof}{{\sc Proof.} \quad}
\newcommand{\proofc}{{\sc Proof} \ }
\newcommand{\be}{\begin{equation} \label}
\newcommand{\ee}{\end{equation}}
\newcommand{\bea}{\begin{eqnarray}\label}
\newcommand{\eea}{\end{eqnarray}}
\newcommand{\bas}{\begin{eqnarray*}}
\newcommand{\eas}{\end{eqnarray*}}
\newcommand{\bit}{\begin{itemize}}
\newcommand{\eit}{\end{itemize}}
\newcommand{\qed}{\hfill$\Box$ \vskip.2cm}
\newcommand{\nn}{\nonumber}
\newcommand{\R}{\mathbb{R}}
\newcommand{\N}{\mathbb{N}}
\newcommand{\pO}{\partial\Omega}

\newcommand{\eps}{\varepsilon}

\newcommand{\wto}{\rightharpoonup}

\newcommand{\io}{\int_\Omega}
\newcommand{\iint}{\int \!\!\! \int}

\newcommand{\abs}{\\[5pt]}

\newcommand{\pos}{{\cal P}}
\newcommand{\tm}{T_{max}}

\newcommand{\Chi}{\chi_{\{u>0\}}}
\begin{document}
%
%
%
%
%
\title{A degenerate fourth-order
parabolic equation \\ modeling Bose-Einstein condensation. \\
Part I: Local existence of solutions} 

\author{
Ansgar J\"ungel\footnote{juengel@tuwien.ac.at}\\
{\small Institute for Analysis and Scientific Computing, Vienna University of Technology,}\\
{\small Wiedner Hauptstra\ss e 8--10, 1040 Wien, Austria} 
\and
Michael Winkler\footnote{michael.winkler@math.uni-paderborn.de}\\	
{\small Institut f\"ur Mathematik, Universit\"at Paderborn,}\\
{\small 33098 Paderborn, Germany} }

\date{}
\maketitle

\begin{abstract}
\noindent
A degenerate fourth-order parabolic equation modeling condensation phenomena
related to Bose-Einstein particles is analyzed. The model is a Fokker-Planck-type
approximation of the Boltzmann-Nordheim equation, only keeping the leading order 
term. It maintains some of the main features of the kinetic model, namely
mass and energy conservation and condensation at zero energy. 
The existence of a local-in-time nonnegative continuous weak solution is proven. 
If the solution is not global, it blows up with respect to the $L^\infty$ norm 
in finite time. 
The proof is based on approximation arguments, interpolation inequalities in
weighted Sobolev spaces, and suitable a priori estimates
for a weighted gradient $L^2$ norm.\abs
\end{abstract}


\mysection{Introduction}

The dynamics of weakly interacting quantum particles like bosons can be described 
by the homogeneous Boltzmann-Nordheim equation for the distribution function $f(x,t)$ 
depending on the energy $x\ge 0$ and time $t>0$ \cite{JPR06},
\begin{equation}\label{1.bne}
  f_t(x_1,t) = \frac{1}{\sqrt{x_1}}\int_{D}S
	\big(f_3 f_4(1+f_1)(1+f_2) - f_1f_2(1+f_3)(1+f_4)\big)dx_3dx_4,
\end{equation}
where $x_2=x_3+x_4-x_1$, $D=\{x_3+x_4>x_1\}$, and the transition rate $S$ 
in the energy space depends on $x_1,\ldots,x_4$. The main feature of this
equation is the existence of finite-time blow-up solutions if the 
initial density is sufficiently dense, modeling the condensation process 
\cite{EsVe12}. The post-nucleation self-similar solution was investigated in detail
by Spohn \cite{Spo10}. Due to the high complexity of the Boltzmann-Nordheim
equation, approximate Fokker-Planck-type equations modeling condensation 
phenomena related to Bose-Einstein particles were studied in the literature.\abs
For instance, if the energy exchange of each collision is small, the
Fokker-Planck approximation of the Nordheim equation in the non-relativistic 
regime leads to the so-called Kompaneets equation \cite{Kom57}. It was 
originally suggested to describe the evolution of a homogeneous plasma
when radiation interacts with matter via Compton scattering. 
Escobedo et al.\ \cite{EHV98} showed that this equation develops singularites
at zero energy.\abs
Another Fokker-Planck model was studied by Kaniadakis and Quarati 
\cite{KaQu93,KaQu94}, proposing a nonlinear correction to the linear
drift term to account for the presence of quantum indistinguishable particles
(bosons and fermions). The model was derived  in \cite{AlTo10}
from a Boltzmann Bose-Einstein model in
the crazing collision limit. Toscani \cite{Tos12} proved that the limit
equation possesses global-in-time solutions if the initial mass is
sufficiently small and the solutions blow up in finite time if the initial mass
is large enough.\abs
A Fokker-Planck-type equation, only containing the superlinear
drift term, was analyzed recently by Carrillo et al.\ \cite{CDT13}. The existence of
a unique measure-valued solution, which concentrates the mass at the origin, 
was proven. Moreover, all mass concentrates in the long-time limit $t\to\infty$.\abs
All these Fokker-Planck equations are of first or second order. 
A {\em higher-order} 
Fokker-Planck approximation of the Boltzmann-Nordheim equation was motivated
by Josserand et al.\ \cite{JPR06}. This model is the subject of this paper.
Assuming that the main contribution to the
collision operator on the right-hand side of (\ref{1.bne}) comes from the
neighborhood of $x\approx x_1\approx x_2\approx x_3\approx x_4$, the integrand of
the collision operator can be expanded to second order, leading to the
fourth-order parabolic equation
\begin{equation}\label{1.fp}
  u_t = x^{-1/2}\Big(x^{13/2}\big(u^{4}(u^{-1})_{xx}-u^2(\log u)_{xx}\big)\Big)_{xx},
	\quad x\in(0,\infty),\ t>0,
\end{equation}
where $u(x,t)$ denotes the energy distribution. 
This approximation maintains some
of the features of the original Boltzmann equation. Indeed, 
assuming no-flux-type boundary
conditions at $x=0$ and $x\to\infty$, this equation conserves the total mass
$N=\int_0^\infty x^{1/2}u dx$ and the kinetic energy $E=\int_0^\infty x^{3/2} udx$.
Furthermore, the entropy
$S = \int_0^\infty((1+u)\log(1+u)-u\log u)x^{1/2}dx$ 
is nondecreasing, and the equilibrium is reached at the Bose-Einstein distribution
$u=(e^{(x-\mu)/T}-1)^{-1}$, where $\mu$ and $T$ are some parameters \cite{JPR06}.\abs
We expect that the local approximation (\ref{1.fp}) contains the relevant 
information on the finite-time collapse of the distribution function. 
For such a study, it is reasonable
to keep only the leading-order cubic term in (\ref{1.fp}). 
Furthermore, we restrict ourselves to the finite energy interval $(0,L)$ for
an arbitrarily large $L>0$ to avoid some technicalities due to infinite domains.
Because of the condensation at energy $x=0$, we expect that the density 
essentially vanishes for large energies which makes Neumann-type boundary conditions
at $x=L$ plausible.\abs
More precisely, in this paper we shall subsequently consider the slightly generalized
problem given by
\be{0}
	\left\{ \begin{array}{rl}
	u_t= x^{-\beta} \Big( x^\alpha u^{n+2} (u^{-1})_{xx} \Big)_{xx},
	& \qquad x\in \Omega, \ t>0, \\[2mm]
	x^\alpha u^{n+2} (u^{-1})_{xx}= \Big( x^\alpha u^{n+2} (u^{-1})_{xx} \Big)_{xx}=0,
	& \qquad x=0, \ t>0, \\[2mm]
	u_x=u_{xxx}=0, & \qquad x=L, \ t>0, \\[2mm]
	u(x,0)=u_0(x), & \qquad x\in \Omega,
	\end{array} \right.
\ee
where $\alpha\ge 0$, $\beta\in\R$, $n>0$, and $\Omega=(0,L)\subset \R$,
with a given nonnegative function $u_0$.\\
The boundary conditions at $x=0$ correspond to those imposed in 
\cite[Formulas (13)-(14)]{JPR06}. In the original equation, we have
$\alpha=13/2$, $\beta=1/2$, and $n=2$. The approximate equation in (\ref{0})
still conserves mass and energy. Moreover, it admits the stationary solutions
$u(x)=x^{-\sigma}$ with $\sigma\in\{0,1,\frac{7}{6},\frac{3}{2}\}$, containing the same
Kolmogorov-Zkharov spectra as the full Boltzmann-Nordheim equation 
\cite[Section 3.3]{JPR06}. This indicates that there is condensation at zero 
energy $x=0$.\abs
{}From a mathematical point of view, 
significant challenges for the analysis stem from the 
fact that the parabolic equation in (\ref{0}) degenerates both at $u=0$ and at $x=0$;
accordingly, the literature does not yet provide any result for this equation,
except for the heuristic study on self-similar solutions in \cite{JPR06}.
It will turn out that this double degeneracy drastically distinguishes the solution behavior in (\ref{0}) 
from that in related well-studied 
degenerate fourth-order parabolic equations such as the thin-film equation
$u_t+(u^n u_{xxx})_x=0$ \cite{BeBeDalPa95,BeGr05,DGG98}. 
Whereas e.g.~the Neumann problem for the latter equation always possesses
a globally defined continuous weak solution which remains bounded \cite{BeFr90,Ber98},
we shall see in the forthcoming paper \cite{JuWi13} that the particular interplay of degeneracies
in (\ref{0}) can enforce solutions to blow up with respect to their spatial norm in
$L^\infty(\Omega)$ within finite time.
More generally, quite various types of higher-order diffusion equations such as e.g.~the 
quantum diffusion or Derrida-Lebowitz-Speer-Spohn equation
\cite{GST09,JuMa08}, equations of epitaxial thin-film growth \cite{Win11}, or also 
some nonlinear sixth-order equations \cite{BJM13,EGK07,PaZa13} have recently attracted considerable interest. To the best of our knowledge, however, such effects of spontaneous singularity formation,
only due to a pure diffusion mechanism without any presence of external forces, have not been detected
in any of these examples.\abs
Against this background, the furthest conceivable outcome of any existence theory can 
only address {\em local} solvability.
The goal of the present work is to establish an essentially optimal result in this direction,
asserting local existence of a continuous weak solution $u$ that conserves mass and that can be extended up to a maximal
existence time $\tm \in (0,\infty]$ at which $\|u(\cdot,t)\|_{L^\infty(\Omega)}$ must blow up
whenever $\tm<\infty$.\abs
Before we state our main result, we introduce some notation. We define
for $\gamma\in\R$ the weighted Sobolev space
$$
  W_\gamma^{1,2}(\Omega) = \big\{v\in W_{\rm loc}^{1,2}(\Omega):
	\|v\|_{L^2(\Omega)}^2+\|x^{\gamma/2}v_x\|_{L^2(\Omega)}<\infty\big\}
$$
with norm $\|v\|_\gamma=(\|v\|_{L^2(\Omega)}^2
+\|x^{\gamma/2}v_x\|_{L^2(\Omega)})^{1/2}$. 
We denote by $\chi_Q$ the characteristic function on the set $Q\subset\R^n$.
The space $C^{4,1}(\bar\Omega\times(0,T))$ consists of all functions
$u$ such that $u_{xxxx}$ and $u_t$ exist and are continuous on 
$\bar\Omega\times(0,T)$. Furthermore, for any (not necessarily open)
subset $Q\subset\R^n$,
$C_0^\infty(Q)$ is the space of all functions such that $\mbox{supp}(f)\subset Q$
is compact.
\begin{defi}\label{defi1}
  Let $n$, $\alpha$, $\beta\in\R$, and $T>0$, and suppose that $u_0\in C^0(\bar\Omega)$ is nonnegative. 
  Then by a {\em continuous weak solution} of (\ref{0}) in
  $\Omega\times (0,T)$ we mean a nonnegative function $u\in C^0(\bar\Omega \times [0,T))$ with the properties
  $u\in C^{4,1}( ((0,L]\times (0,T)) \cap \{u>0\})$ as well as
  \be{0r}
	\Chi x^\alpha u^n u_{xx} \in L^1_{loc}(\bar\Omega \times [0,T))
	\qquad \mbox{and} \qquad
	\Chi x^\alpha u^{n-1} u_x^2 \in L^1_{loc}(\bar\Omega \times [0,T)),
  \ee
  for which $u(\cdot,t)$ is differentiable with respect to $x$ at $x=L$ for a.e.~$t\in (0,T)$ with
  \be{0b}
	u_x(L,t)=0	\qquad \mbox{for a.e.~$t\in (0,T)$},
  \ee
  and which satisfies the integral identity
  \be{0w}
	-\int_0^T \io x^\beta u\phi_t dxdt - \io x^\beta u_0 \phi(\cdot,0)dx
	= \int_0^T \io \chi_{\{u>0\}} [-x^\alpha u^n u_{xx} + 2x^\alpha u^{n-1} u_x^2 ] \phi_{xx}dxdt
  \ee
  for all $\phi \in C_0^\infty(\bar\Omega \times [0,T))$ fulfilling $\phi_x(L,t)=0$ for all $t\in (0,T)$.
\end{defi}
Note that if $u$ is a positive classical solution in the sense of this definition
and $\alpha>1$, then partial integration in (\ref{0w}) shows that
$u$ satisfies the boundary conditions in (\ref{0}).
Our main result reads as follows.
\begin{theo}[Local existence of solutions]\label{theo_final}
  Let $n\in(n^*,3)$, where $n^*=1.5361\ldots$ is the unique positive root of 
  the polynomial $n\mapsto n^3+5n^2+16n-40$. Let $\alpha>3$ and
  $\beta\in(-1,\alpha-4)$. Then for any $\gamma\in(5-\alpha+\beta,1)$ and
  each nonnegative function $u_0\in W_\gamma^{1,2}(\Omega)$, there exists
  $\tm\in(0,\infty]$ such that (\ref{0}) possesses 
  a continuous weak solution $u\in L_{loc}^\infty([0,\tm);W^{1,2}_\gamma(\Omega))$. 
  Furthermore, 
  \begin{equation}\label{extend}
 	\mbox{if }\tm<\infty\mbox{ then }
  	\limsup_{t\to\tm}\|u(\cdot,t)\|_{L^\infty(\Omega)}=\infty,
  \end{equation}
  and the solution conserves the mass in the sense that
  \bas
	\int_\Omega x^\beta u(x,t)dx = \int_\Omega x^\beta u_0(x)dx
	\quad\mbox{for a.e. }t\in(0,\tm).
  \eas
\end{theo}
Note that the physical values $\alpha=\frac{13}{2}$, $\beta=\frac{1}{2}$, and $n=2$ are admissible
choices in the theorem.\abs
A cornerstone in our analysis will consist in establishing an a priori estimate of the form
\bea{formal}
	\frac{d}{dt} \io x^\gamma u_x^2dx
	&+& c \io x^{\alpha-\beta+\gamma} u^n u_{xxx}^2 dx
	+ c \io x^{\alpha-\beta+\gamma} u^{n-2} u_x^2 u_{xx}^2 dx
	+ c \io x^{\alpha-\beta+\gamma} u^{n-4} u_x^6 dx \nn\\
	&+& c \io x^{\alpha-\beta+\gamma-2} u^n u_{xx}^2 dx
	+ c \io x^{\alpha-\beta+\gamma-2} u^{n-2} u_x^4 dx \nn\\
	&\le& C + C \bigg( \io x^\gamma u_x^2 dx\bigg)^\frac{n+2}{2}
\eea
for appropriate $c>0$ and $C>0$,	
which can formally be derived from (\ref{0}) under the restrictions for $\alpha$,
$\beta$, $\gamma$, and $n$ made in Theorem \ref{theo_final}.
Upon integration, (\ref{formal}) will imply appropriate weighted integral estimates for $u$ and its derivatives
on small time intervals, inter alia the inequality
\be{formal2}
	\io x^\gamma u_x^2(x,t)dx \le \tilde C \qquad \mbox{for all } t\in (0,T)
\ee
for some $\tilde C>0$ and appropriately small $T>0$.\\
A rigorous variant of (\ref{formal}) is shown in Lemmas \ref{lem7} and 
\ref{lem13}. In view of the degeneracies in (\ref{0}), our analysis will
rely on a suitable regularization. To achieve this,
we shall replace $x^{-\beta}$ and $x^\alpha$ by $(x+\eps)^{-\beta}$ and 
$g_\eps(x)$, respectively, where $\eps>0$,
$g_\eps$ is positive in $\Omega$, and $g_{\eps,x}$ vanishes on the boundary.
The latter condition ensures that the approximate flux 
$J=-g_\eps(x)(-u^n u_{xx}+2u^{n-1}u_x^2)$ vanishes on the boundary as well.
We emphasize that unlike typical approaches in related equations such as the thin-film equation,
our regularized problems are still degenerate at $u=0$. To circumvent obstacles stemming from this,
we shall first consider stricly positive initial data only; 
however, this will require additional efforts in ruling out that the local-in-time approximate solutions thereby obtained
do not approach this critical level $u=0$ within finite time 
(see Lemma \ref{lem12}).\\
The limit process $\eps\to 0$ will then be carried out on the basis of a spatio-temporal H\"older estimate for
the approximate solutions, which thanks to the fact that $\gamma<1$ can be derived from (\ref{formal2})
along with the adaptation of a well-known argument from parabolic theory, which turns this into an appropriate H\"older
estimate with respect to time (Lemma \ref{lem14}).\abs
The paper is organized as follows. In Section \ref{sec.approx}, we introduce
the family of approximate problems. Interpolation inequalities in weighted spaces,
which are needed for the existence analysis, are shown in Section \ref{sec.inter}.
The proof of the a priori estimates (Lemmas \ref{lem7} and \ref{lem13})
is the subject of Section \ref{sec.est}. Then Section \ref{sec.locex} is concerned
with the local existence for the approximate problems and the absence of dead core
formation. A H\"older estimate for the approximate solutions is derived in
Section \ref{sec.hoelder}. Finally, the proof of Theorem \ref{theo_final} is presented in 
Section \ref{sec.proof}.


\section{A family of approximate problems}\label{sec.approx}
We formulate a family of approximate problems in which the
singularity at $x=0$ is removed but the boundary conditions in (\ref{0})
hold at $x=L$ and $x=0$.
To this end, we let $\eps_0=\min\{1,\sqrt{L/2}\}$, and for $\eps\in(0,\eps_0)$,
we choose $\zeta_\eps\in C_0^\infty(\Omega)$ satisfying $0\le\zeta_\eps\le 1$
and $\zeta_\eps(y)=1$ for $y\in(\eps^2,L-\eps^2)$. Furthermore, we set
$$
  z_\eps(x) = \eps + \int_0^x \zeta_\eps(y)dy, \quad x\in[0,L].
$$
Then the function $z_\eps$ belongs to $C^\infty([0,L])$, 
$z_\eps(x)\ge \eps$ for all $x\in[0,L]$, and it satisfies homogeneous
Neumann boundary conditions, $z_{\eps,x}(0)=z_{\eps,x}(L)=0$.
Then $g_\eps:=z_\eps^\alpha$ belongs to $C^\infty([0,L])$ and satisfies $g_\eps\ge \eps^\alpha$ on $[0,L]$ and
$g_{\eps,x}(0)=g_{\eps,x}(L)=0$. Further pointwise estimates for $g_\eps$ are 
summarized in the following lemma.

\begin{lem}[Properties of $g_\eps$]\label{lem1}
  Let $\alpha>0$. Then the following properties hold:\\
  (i) \ There exists a positive decreasing function $\Lambda:[0,\eps_0)\to(0,1)$
  such that $\inf_{(0,\eps_0)}\Lambda>0$, $\Lambda(0)=1$,
  and for all $\eps\in(0,\eps_0)$,
  \bas
	\Lambda(\eps)(x+\eps)^\alpha \le g_\eps(x)\le (x+\eps)^\alpha, \quad x\in[0,L].  
  \eas
  (ii) \ There exists $c>0$ such that for all $\eps\in(0,\eps_0)$,
  \bas
	0 \le g_{\eps x}(x)\le c(x+\eps)^{\alpha-1}, \quad x\in[0,L].
  \eas
  (iii) \ There exists $c>0$ such that for all $\eps\in(0,\eps_0)$,
  \bas
	\frac{g_{\eps x}(x)^2}{g_\eps(x)} \le c(x+\eps)^{\alpha-2}, \quad x\in[0,L].
  \eas
\end{lem}
\proof
(i) \ Since $\zeta_\eps\le 1$, we have $z_\eps(x)\le \eps+x$ for $x\in[0,L]$. 
This yields the second inequality, $g_\eps(x)=z_\eps(x)^\alpha\le (x+\eps)^\alpha$.
To prove the first one, 
we divide $[0,L]$ into three subintervals. First, for $x\in[0,\eps^2]$
the property $z_\eps(x)\ge \eps$ yields
$$
  \frac{z_\eps(x)}{x+\eps} \ge \frac{\eps}{x+\eps}\ge \frac{\eps}{\eps^2+\eps}
	= \frac{1}{1+\eps}.
$$
Next, if $x\in(\eps^2,L-\eps^2)$ then $\zeta_\eps(x)=1$, whence using that
$z_\eps(\eps^2)\ge\eps$, we obtain
\bas
  \frac{z_\eps(x)}{x+\eps}
	&=& \frac{1}{x+\eps}\left(z_\eps(x)+\int_{\eps^2}^x\zeta_\eps(y)dy\right)
	\ge \frac{\eps + (x-\eps^2)}{x+\eps} 
	\ge 1-\frac{\eps^2}{\eps^2+\eps}
	=\frac{1}{1+\eps}.
\eas
We finally consider the case $x\in[L-\eps^2,L]$, in which because of the nonnegativity of $z_\eps$ and the fact that
$\zeta_\eps=1$ on $[\eps^2,L-\eps^2]$, we infer that
$$
  \frac{z_\eps(x)}{x+\eps}
	\ge \frac{z_\eps(L-\eps^2)}{x+\eps}
	= \frac{\eps+(L-2\eps^2)}{x+\eps}
	\ge \frac{L+\eps-2\eps^2}{L+\eps} = 1-\frac{2\eps^2}{L+\eps}.
$$
The claim hence follows by defining $\Lambda(\eps)=\min\{1/(1+\eps),1-2\eps^2/(L+\eps)\}$.\abs
(ii) \ As $0\le z_{\eps,x}\le 1$, we have $g_{\eps,x}=\alpha z_\eps^{\alpha-1}z_{\eps,x}
\le \alpha z_\eps^{\alpha-1}$ in $[0,L]$. Thus, (i) implies (ii).\abs
(iii) \ This follows directly from (i) and (ii).
\qed
With the above choices of $\eps_0>0$ and $g_\eps$, we proceed to regularize the original problem appropriately. 
The idea is to replace in the first equation in (\ref{0}),
rewritten in the form $u_t=x^{-\beta}(-x^\alpha u^n u_{xx}+2x^\alpha u^{n-1}u_x^2)_{xx}$,
the coefficients $x^{-\beta}$ and $x^\alpha$ by $(x+\eps)^{-\beta}$ and
$g_\eps(x)$, respectively. Accordingly, for $\eps\in (0,\eps_0)$, we shall consider the approximate problem
\be{0eps}
	\left\{ \begin{array}{l}
	u_t=\frac{1}{(x+\eps)^\beta} \cdot \Big\{ -g_\eps(x) u^n u_{xx} + 2g_\eps(x) u^{n-1} u_{xx} \Big\}_{xx},
	\qquad x\in \Omega, \ t>0, \\[1mm]
	u_x=u_{xxx}=0, \qquad x\in \partial\Omega, \ t>0.
	\end{array} \right.
\ee
The boundary behavior of $g_\eps$ guarantees that the flux
\begin{equation}\label{J}
  J(x,t) = -g_\eps(x)u^n u_{xx} + 2g_\eps(x)u^{n-1}u_x^2
\end{equation}
vanishes on $\pO=\{0,L\}$. 
This results upon expanding $J_x$ according to
\bea{Jx}
	J_x &=& -g_\eps(x) u^n u_{xxx} + (4-n) g_\eps(x) u^{n-1} u_x u_{xx} 
	+ 2(n-1) g_\eps(x) u^{n-2} u_x^3 \nn\\
	& &- g_{\eps,x}(x) u^n u_{xx} + 2g_{\eps,x}(x) u^{n-1} u_x^2,
\eea
and evaluating this expression on $\pO$:
\begin{lem}[Boundary flux vanishes]\label{lem_J}
  Let $n>0$, $\alpha>0$, $\beta\in\R$, $T>0$, and $\eps\in(0,\eps_0)$,
  and let $u\in C^{4,1}(\overline\Omega\times(0,T))$
  be a positive classical solution of (\ref{0eps}). Then
  $J_x(x,t)=0$ for all $x\in\pO$ and $t\in(0,T)$, where $J$ is defined in (\ref{J}).
\end{lem}
\proof
  The statement is a consequence of (\ref{Jx}) and the identities $u_x=u_{xxx}=g_{\eps,x}=0$ on $\pO$.
\qed
The above choice of boundary conditions ensures that the total mass is preserved.
\begin{lem}[Conservation of total mass]\label{lem8}
  Under the assumptions of Lemma \ref{lem_J}, we have 
  \bas
	\frac{d}{dt}\int_\Omega (x+\eps)^\beta u(x,t)dx = 0
	\qquad \mbox{for all } t\in (0,T).
  \eas
\end{lem}
\proof
  The claim immediately results by integrating (\ref{0eps}) over $\Omega$ and using that $J_x=0$ on $\pO$.
\qed
\mysection{Some interpolation inequalities}\label{sec.inter}
As a preparation for our subsequent analysis, let us collect some interpolation inequalities in weighted spaces. 
The first of these reads as follows.
\begin{lem}\label{lem_inter}
  Let $n \in \R \setminus \{-1,1\}$, $\alpha\in\R$, $\beta\in\R$, and $\gamma\in\R$. 
  Then for any $\eta>0$, one can find $C(\eta)>0$ such that for all positive functions $u\in C^2(\bar\Omega)$
  satisfying $u_x=0$ on $\pO$, we have
  \bea{inter1}
	& & \hspace*{-40mm}
	\io (x+\eps)^{\alpha-\beta+\gamma-2}u^{n} u_{xx}^2 dx
	+ \int_\Omega(x+\eps)^{\alpha-\beta+\gamma-2}u^{n-2} u_{x}^4dx
	+ \int_\Omega(x+\eps)^{\alpha-\beta+\gamma-4}u^n u_x^2dx \nn\\
	&\le& 
  	\eta\int_\Omega(x+\eps)^{\alpha-\beta+\gamma} u^n u_{xxx}^2dx
	+ \eta\int_\Omega(x+\eps)^{\alpha-\beta+\gamma} u^{n-2}u_x^2 u_{xx}^2 dx \nn\\
	& & {}+ C(\eta) \io (x+\eps)^{\alpha-\beta+\gamma-6} u^{n+2}dx
  \eea
  for all $\eps>0$.
\end{lem}
The proof of Lemma \ref{lem_inter} will be achieved in a series of steps to be presented separately in Lemmas \ref{lem2}-\ref{lem4}.
We first estimate the last integral on the left-hand side of (\ref{inter1}) 
by a sum involving a small portion of the first term in (\ref{inter1}).
\begin{lem}\label{lem2}
  Let $n \in \R \setminus \{-1\}$ and $\alpha$, $\beta$, and $\gamma$ be 
	arbitrary real numbers. 
  Then for all $\eta>0$, there exists $C(\eta)>0$ such that whenever $\eps>0$ and 
  $u\in C^2(\bar\Omega)$ is positive with $u_x=0$ on $\pO$,
  the inequality
  \be{2.1}
	\io (x+\eps)^{\alpha-\beta+\gamma-4} u^n u_x^2dx
	\le \eta \io (x+\eps)^{\alpha-\beta+\gamma-2} u^n u_{xx}^2dx
	+ C(\eta) \io (x+\eps)^{\alpha-\beta+\gamma-6} u^{n+2}dx
  \ee
  holds.
\end{lem}
\proof
  Using $u_x=0$ on $\pO$, we may integrate by parts and use Young's inequality 
  to find that
  \bas
	\Gamma&:=&\io (x+\eps)^{\alpha-\beta+\gamma-4} u^n u_x^2dx \\
	&=& -\frac{1}{n+1} \io (x+\eps)^{\alpha-\beta+\gamma-4} u^{n+1} u_{xx}dx
	- \frac{\alpha-\beta+\gamma-4}{n+1} \io (x+\eps)^{\alpha-\beta+\gamma-5} u^{n+1} u_xdx \nn\\
	&\le& \frac{\eta}{2} \io (x+\eps)^{\alpha-\beta+\gamma-2} u^n u_{xx}^2dx
	+ \frac{1}{2(n+1)^2 \eta} \io (x+\eps)^{\alpha-\beta+\gamma-6} u^{n+2}dx \nn\\
	& & {}+ \frac{1}{2} \Gamma 
	+ \frac{(\alpha-\beta+\gamma-4)^2}{2(n+1)^2} \io (x+\eps)^{\alpha-\beta+\gamma-6} u^{n+2}dx.
  \eas
  Rearranging yields (\ref{2.1}).
\qed
Using the above preparation, we can control the first term in (\ref{inter1}) as desired:
\begin{lem}\label{lem3}
  Let $n\in \R \setminus \{-1\}$ and $\alpha$, $\beta$, $\gamma\in\R$. 
	Then for all $\eta>0$, one can find
  $C(\eta)>0$ with the property that any positive function $u\in C^3(\bar\Omega)$ with  
	$u_x=0$ on $\pO$ satisfies
  \bea{3.1}
	\io (x+\eps)^{\alpha-\beta+\gamma-2} u^n u_{xx}^2dx
	&\le& \eta \io (x+\eps)^{\alpha-\beta+\gamma} u^n u_{xxx}^2 dx
	+ \eta \io (x+\eps)^{\alpha-\beta+\gamma} u^{n-2} u_x^2 u_{xx}^2dx \nn\\
	& & {}+ C(\eta) \io (x+\eps)^{\alpha-\beta+\gamma-6} u^{n+2}dx
  \eea
  for each $\eps>0$.
\end{lem}
\proof
  Since $u_x=0$ on $\pO$, an integration by parts shows that
  \bea{3.2}
	\Gamma&:=& \io (x+\eps)^{\alpha-\beta+\gamma-2} u^n u_{xx}^2dx
	= - \io (x+\eps)^{\alpha-\beta+\gamma-2} u^n u_x u_{xxx} dx \\
	&&{}- n \io (x+\eps)^{\alpha-\beta+\gamma-2} u^{n-1} u_x^2 u_{xx}dx 
	- (\alpha-\beta+\gamma-2) \io (x+\eps)^{\alpha-\beta+\gamma-3} u^n u_x u_{xx}dx,
	\nn
  \eea
  where by Young's inequality we find that
  \bas
	- \io (x+\eps)^{\alpha-\beta+\gamma-2} u^n u_x u_{xxx} dx
	\le \frac{\eta}{2} \io (x+\eps)^{\alpha-\beta+\gamma} u^n u_{xxx}^2dx
	+ c_1 \io (x+\eps)^{\alpha-\beta+\gamma-4} u^n u_x^2dx
  \eas
  and
  \bas
	- n \io (x+\eps)^{\alpha-\beta+\gamma-2} u^{n-1} u_x^2 u_{xx} dx
	&\le& \frac{\eta}{4} \io (x+\eps)^{\alpha-\beta+\gamma} u^{n-2} u_x^2 u_{xx}^2dx \\
	&&{}+ c_2 \io (x+\eps)^{\alpha-\beta+\gamma-4} u^n u_x^2dx
  \eas
  as well as
  \bas
	-(\alpha-\beta+\gamma-2) \io (x+\eps)^{\alpha-\beta+\gamma-3} u^n u_x u_{xx}dx
	&\le& \frac{\eta}{4} \io (x+\eps)^{\alpha-\beta+\gamma} u^{n-2} u_x^2 u_{xx}^2dx \\
	& & {}+ c_3 \io (x+\eps)^{\alpha-\beta+\gamma-6} u^{n+2}dx
  \eas
  with $c_1:=\frac{1}{2\eta}, c_2:=\frac{n^2}{\eta}$ and $c_3:=\frac{(\alpha-\beta+\gamma-2)^2}{\eta}$.
  Since Lemma \ref{lem2} provides $c_4>0$ such that
  \bas
	(c_1+c_2) \io (x+\eps)^{\alpha-\beta+\gamma-4} u^n u_x^2 dx
	\le \frac{1}{2} \Gamma + c_4 \io (x+\eps)^{\alpha-\beta+\gamma-6} u^{n+2}dx,
  \eas
  (\ref{3.2}) thereby proves (\ref{3.1}).
\qed
%
%
%
%
%
%
%
%
%
%
%
%
Now the latter allows us to also estimate the second term in (\ref{inter1}) in the claimed manner.
\begin{lem}\label{lem4}
  Let $n\in\R\setminus\{-1,1\}$ and $\alpha$, $\beta$, $\gamma \in \R$. Then for 
	all $\eta>0$, we can pick $C(\eta)>0$ such that if 
  $u\in C^3(\bar\Omega)$ is positive and satisfies $u_x=0$ on $\pO$, 
	then for all $\eps>0$, we have
  \bea{4.1}
	\io (x+\eps)^{\alpha-\beta+\gamma-2} u^{n-2} u_x^4dx
	&\le& \eta \io (x+\eps)^{\alpha-\beta+\gamma} u^n u_{xxx}^2dx
	+ \eta \io (x+\eps)^{\alpha-\beta+\gamma} u^{n-2} u_x^2 u_{xx}^2dx \nn\\
	& & {}+ C(\eta) \io (x+\eps)^{\alpha-\beta+\gamma-6} u^{n+2}dx.
  \eea
\end{lem}
\proof
  Once more integrating by parts and using Young's inequality, we see that
  \bas
	\Gamma &:=& \io (x+\eps)^{\alpha-\beta+\gamma-2} u^{n-2} u_x^4dx \\
	&=& - \frac{3}{n-1} \io (x+\eps)^{\alpha-\beta+\gamma-2} u^{n-1} u_x^2 u_{xx}dx 
	- \frac{\alpha-\beta+\gamma-2}{n-1} \io (x+\eps)^{\alpha-\beta+\gamma-3} u^{n-1} u_x^3dx \nn\\
	&\le& \frac{1}{4} \Gamma
	+ c_1 \io (x+\eps)^{\alpha-\beta+\gamma-2} u^n u_{xx}^2dx
	+ \frac{1}{4} \Gamma
	+ c_2 \io (x+\eps)^{\alpha-\beta+\gamma-6} u^{n+2}dx
  \eas
  with $c_1:=\frac{9}{(n-1)^2}$ and $c_2:=\frac{(\alpha-\beta+\gamma-4)^2}{(n-1)^2}$. Thus,
  \bas
	\Gamma \le 2c_1 \io (x+\eps)^{\alpha-\beta+\gamma-2} u^n u_{xx}^2dx
	+ 2c_2 \io (x+\eps)^{\alpha-\beta+\gamma-6} u^{n+2}dx,
  \eas
  whence invoking Lemma \ref{lem3}, we readily arrive at (\ref{4.1}).
\qed
\proofc of Lemma \ref{lem_inter}. \quad
  We only need to combine Lemmas \ref{lem2}, \ref{lem3}, and \ref{lem4}.
\qed
The following inequality is closely related to those used in the context of the 
thin-film equation $u_t+(u^n u_{xxx})_x=0$ \cite{Be96}.
\begin{lem}\label{lem6}
  Let $n\in\R \setminus \{3\}$ and $\alpha$, $\beta$, $\gamma \in\R$. 
	Then for all $\eta\in (0,1)$ and
  any positive $u\in C^2(\bar\Omega)$ fulfilling $u_x=0$ on $\pO$, the inequality
  \bea{6.1}
	\io (x+\eps)^{\alpha-\beta+\gamma} u^{n-4} u_x^6dx
	&\le& \frac{25}{(1-\eta)(n-3)^2} \io (x+\eps)^{\alpha-\beta+\gamma} u^{n-2} u_x^2 u_{xx}^2dx \nn\\
	& & {}+ \frac{(\alpha-\beta+\gamma)^2}{\eta(1-\eta)(n-3)^2} \io (x+\eps)^{\alpha-\beta+\gamma-2} u^{n-2} u_x^4dx
  \eea
  is valid for all $\eps>0$.
\end{lem}
\proof
  We integrate by parts using $u_x=0$ on $\pO$ and apply Young's inequality to obtain the estimate
  \bas
	\Gamma&:=&
	\io (x+\eps)^{\alpha-\beta+\gamma} u^{n-4} u_x^6dx \\
	&=& - \frac{5}{n-3} \io (x+\eps)^{\alpha-\beta+\gamma} u^{n-3} u_x^4 u_{xx}dx
  - \frac{\alpha-\beta+\gamma}{n-3} \io (x+\eps)^{\alpha-\beta+\gamma-1} u^{n-3} u_x^5dx \\
	&\le& \frac{1}{2} \Gamma
	+ \frac{1}{2} \cdot \frac{25}{(n-3)^2} \io (x+\eps)^{\alpha-\beta+\gamma} u^{n-2} u_x^2 u_{xx}^2dx \\
	&& {}+ \frac{\eta}{2} \Gamma
	+ \frac{1}{2\eta} \cdot \frac{(\alpha-\beta+\gamma)^2}{(n-3)^2} \io (x+\eps)^{\alpha-\beta+\gamma-2} u^{n-2} u_x^4dx,
  \eas
  which can readily be checked to be equivalent to (\ref{6.1}).
\qed
The following two lemmas are concerned with estimates on the H\"older and
$L^\infty$ norms of functions in $W^{1,2}_{loc}(\Omega)$.
\begin{lem}\label{lem31}
  Let $\gamma\in(-\infty,1)$. Then there exists $c(\gamma)>0$
  such that for any $\eps\in[0,1)$ and any $u\in W^{1,2}_{loc}(\Omega)$,
  \bas
  	|u(x_2)-u(x_1)|\le c\left(\int_\Omega(x+\eps)^\gamma u_x^2 dx\right)^\frac{1}{2}
	|x_2-x_1|^\theta \quad\mbox{for all }x_1,x_2\in\Omega,
  \eas
  where $\theta:=\min\{\frac{1}{2}, \frac{1-\gamma}{2}\}$.
\end{lem}
\proof
  Let $0<x_1<x_2<L$ and suppose first that $\gamma\in[0,1)$. Then by the Cauchy-Schwarz inequality,
  \bas
	|u(x_2)-u(x_1)| = \left|\int_{x_1}^{x_2}u_x(x)dx\right|
	\le \left(\int_\Omega(x+\eps)^\gamma u_x^2 dx\right)^\frac{1}{2}
	\left(\int_{x_1}^{x_2}(x+\eps)^{-\gamma}dx\right)^\frac{1}{2}.
  \eas
  Employing the H\"older continuity of $x\mapsto x^{1-\gamma}$, we obtain
  \bas
	\int_{x_1}^{x_2}(x+\eps)^{-\gamma}dx
	= \frac{1}{1-\gamma}\big((x_2+\eps)^{1-\gamma}-(x_1+\eps)^{1-\gamma}\big) 
	\le \frac{c_1}{1-\gamma}|x_2-x_1|^{1-\gamma}.
  \eas
  The result thus follows with $c=(c_1/(1-\gamma))^\frac{1}{2}$.\\
  If $\gamma\in(-\infty,0)$, we replace $\gamma$ by $-\gamma$ in the above arguments and use the Lipschitz 
  continuity of $x\mapsto x^{1+|\gamma|}$.
\qed
\begin{lem}\label{lem32}
  Let $\gamma\in(-\infty,1)$ and $\beta\in\R$. Then there exists $c=c(\beta,L)>0$
  such that for all $\eps\in[0,1)$ and any $u\in W^{1,2}_{loc}(\Omega)$,
  \bas
	\|u\|_{L^\infty(\Omega)} \le c\left(\int_\Omega (x+\eps)^\beta|u|dx
	+ \left(\int_\Omega(x+\eps)^\gamma u_x^2 dx\right)^{1/2}\right).
  \eas
\end{lem}
\proof
  Assuming that $B=\int_\Omega(x+\eps)^\beta|u|dx$ is finite, we see that there exists
  $x_0\in(\frac{L}{2},L)$ such that $(x_0+\eps)^\beta|u(x_0)|\le \frac{2B}{L}$, for otherwise
  the inequality $B\ge \int_\frac{L}{2}^L(x+\eps)^\beta|u|dx>\frac{L}{2} \cdot \frac{2B}{L}=B$ gives
  a contradiction. Since $\frac{L}{2}\le x_0+\eps\le L+1$, we infer that
  \bas
  	|u(x_0)| \le c_1\int_\Omega(x+\eps)^\beta|u|dx,
  \eas
  where $c_1=\frac{2}{L} \cdot \max\{(\frac{L}{2})^{-\beta},(L+1)^{-\beta}\}$.
  The conclusion thus follows from Lemma \ref{lem31}.
\qed
\mysection{A differential inequality for $\io x^\gamma u_x^2$}\label{sec.est}
A key role in our analysis will be played by the following a priori estimate for the functional
$y(t):=\int_\Omega(x+\eps)^\beta u_x^2dx$ in terms of a weighted norm of $u$ in $L^{n+2}(\Omega)$.
In Lemma \ref{lem13} below,
we shall turn this into an autonomous differential equation for $y(t)$, which will be essential for our local existence proof.
\begin{lem}[A priori estimate in terms of a weighted $L^{n+2}$ norm]\label{lem7}
  Let $n_\star=1.5361\ldots$ be the unique positive root of $n\mapsto P(n):=n^3+5n^2+16n-40$,
  and let $n\in(n_\star,3)$, $\alpha>0$, $\beta\in\R$, and $\gamma\in\R$.
  Then there exist $\eps_\star\in (0,\eps_0)$, 
  $c>0$, and $K>0$ such that if for some $T>0$ and $\eps \in (0,\eps_\star)$,
  $u\in C^{4,1}(\bar\Omega \times (0,T))$ is a positive classical solution to (\ref{0eps}), then
  \bea{7.1}
	\frac{d}{dt} \io (x+\eps)^\gamma u_x^2 (x,t)dx
	&+& c \io (x+\eps)^{\alpha-\beta+\gamma} u^n u_{xxx}^2dx
	+ c \io (x+\eps)^{\alpha-\beta+\gamma} u^{n-2} u_x^2 u_{xx}^2dx \nn\\
	&+& c \io (x+\eps)^{\alpha-\beta+\gamma-2} u^n u_{xx}^2dx
	+ c \io (x+\eps)^{\alpha-\beta+\gamma-2} u^{n-2} u_x^4dx \nn\\
	&\le& K \io (x+\eps)^{\alpha-\beta+\gamma-6} u^{n+2}dx
	\qquad \mbox{for all } t\in (0,T).
  \eea
\end{lem}
\proof
  With the notation (\ref{J}), we can write the first equation in (\ref{0eps}) as
  $u_t=(x+\eps)^{-\beta}J_{xx}$. Since $u_x=J_x=0$ on $\pO$ by Lemma \ref{lem_J}, an integration by parts gives
  \bas
  	\frac{1}{2} \, \frac{d}{dt}\int_\Omega(x+\eps)^\gamma u_x^2 dx
	&=& -\int_\Omega((x+\eps)^\gamma u_x)_x u_t dx
	= -\int_\Omega((x+\eps)^\gamma u_x)_x (x+\eps)^{-\beta}J_{xx} dx \\
	&=& \int_\Omega\Big[(x+\eps)^{-\beta+\gamma}u_{xx} 
	+ \gamma(x+\eps)^{-\beta+\gamma-1}u_x\Big]_x J_x dx
  \eas
  for all $t\in (0,T)$.
  Computing
  \bas
	& & \hspace*{-20mm}
	\Big[ (x+\eps)^{-\beta+\gamma} u_{xx} + \gamma (x+\eps)^{-\beta+\gamma-1} u_x \Big]_x \\[2mm]
	&=& (x+\eps)^{-\beta+\gamma} u_{xxx} + (2\gamma-\beta) (x+\eps)^{-\beta+\gamma-1} u_{xx}
	+ \gamma(\gamma-\beta-1) (x+\eps)^{-\beta+\gamma-2} u_x
  \eas
  and expanding $J_x$ by means of (\ref{Jx}), we thus obtain the identity
  \bea{7.2}
	\frac{1}{2} \frac{d}{dt} \io (x+\eps)^\gamma u_x^2dx
	&=& - \io (x+\eps)^{-\beta+\gamma} g_\eps(x) u^n u_{xxx}^2dx \nn\\
	& & {}+ (4-n) \io (x+\eps)^{-\beta+\gamma} g_\eps(x) u^{n-1} u_x u_{xx} u_{xxx}dx \nn\\
	& & {}+ 2(n-1) \io (x+\eps)^{-\beta+\gamma} g_\eps(x) u^{n-2} u_x^3 u_{xxx}dx \nn\\
	& & {}- \io (x+\eps)^{-\beta+\gamma} g_{\eps x}(x) u^n u_{xx} u_{xxx}dx \nn\\
	& & {}+ 2 \io (x+\eps)^{-\beta+\gamma} g_{\eps x}(x) u^{n-1} u_x^2 u_{xxx}dx \nn\\
	& & {}- (2\gamma-\beta) \io (x+\eps)^{-\beta+\gamma-1} g_\eps(x) u^n u_{xx} u_{xxx} dx\nn\\
	& &{}+(2\gamma-\beta)(4-n) \io (x+\eps)^{-\beta+\gamma-1} g_\eps(x) u^{n-1} u_x u_{xx}^2dx \nn\\
	& & {}+ 2(2\gamma-\beta)(n-1) \io (x+\eps)^{-\beta+\gamma-1} g_\eps(x) u^{n-2} u_x^3 u_{xx}dx \nn\\
	& & {}- (2\gamma-\beta) \io (x+\eps)^{-\beta+\gamma-1} g_{\eps x}(x) u^n u_{xx}^2dx \nn\\
	& & {}+2(2\gamma-\beta) \io (x+\eps)^{-\beta+\gamma-1} g_{\eps x}(x) u^{n-1} u_x^2 u_{xx}dx \nn\\
	& & {}-\gamma(\gamma-\beta-1) \io (x+\eps)^{-\beta+\gamma-2} g_\eps(x) u^n u_x u_{xxx}dx \nn\\
	& & {}+\gamma(\gamma-\beta-1)(4-n) \io (x+\eps)^{-\beta+\gamma-2} g_\eps(x) u^{n-1} u_x^2 u_{xx}dx \nn\\
	& & {}+2\gamma(\gamma-\beta-1)(n-1) \io (x+\eps)^{-\beta+\gamma-2} g_\eps(x) u^{n-2} u_x^4dx \nn\\
	& & {}-\gamma(\gamma-\beta-1) \io (x+\eps)^{-\beta+\gamma-2} g_{\eps x}(x) u^n u_x u_{xx}dx \nn\\
	& & {}+2\gamma(\gamma-\beta-1) \io (x+\eps)^{-\beta+\gamma-2} g_{\eps x}(x) u^{n-1} u_x^3dx \nn\\[2mm]
	&=:& I_1+ \cdots + I_{15}
	\qquad \mbox{for all } t\in (0,T).
  \eea
  Our goal is to adequately apply the interpolation inequalities in Lemma \ref{lem_inter} and Lemma \ref{lem6}
  and to identify those integrals which absorb the $O(\eta)$
  contributions in (\ref{inter1}) and (\ref{6.1}) such that finally only a possibly large multiple of the
  integral over $(x+\eps)^{\alpha-\beta+\gamma-6}u^{n+2}$ remains.\\
  To achieve this, we observe that the integral $I_1$ is nonpositive and thus can be used to absorb positive 
  contributions.
  Apart from this, the only absorptive contribution to be used in the sequel will result from $I_3$, which
  we therefore rearrange first: Namely, by two further integrations by parts, 
  once more relying on the fact that $u_x=0$ on $\pO$, this term can be rewritten according to
  \bea{7.3}
	I_3 &=& -6(n-1) \io (x+\eps)^{-\beta+\gamma} g_\eps(x) u^{n-2} u_x^2 u_{xx}^2dx \nn\\
	& & {}- 2(n-1)(n-2) \io (x+\eps)^{-\beta+\gamma} g_\eps(x) u^{n-3} u_x^4 u_{xx}dx \nn\\
	& & {}-2(\gamma-\beta)(n-1) \io (x+\eps)^{-\beta+\gamma-1} g_\eps(x) u^{n-2} u_x^3 u_{xx}dx \nn\\
	& &{}-2(n-1) \io (x+\eps)^{-\beta+\gamma} g_{\eps x}(x) u^{n-2} u_x^3 u_{xx}dx \nn\\
	&=& -6(n-1) \io (x+\eps)^{-\beta+\gamma} g_\eps(x) u^{n-2} u_x^2 u_{xx}^2dx \nn\\
	& & {}-\frac{2}{5} (n-1)(n-2)(3-n) \io (x+\eps)^{-\beta+\gamma} g_\eps(x) u^{n-4} u_x^6dx \nn\\
	& & {}+ \frac{2}{5}(\gamma-\beta)(n-1)(n-2) \io (x+\eps)^{-\beta+\gamma-1} g_\eps(x) u^{n-3} u_x^5dx \nn\\ 
	& & {}+ \frac{2}{5} (n-1)(n-2) \io (x+\eps)^{-\beta+\gamma} g_{\eps x}(x) u^{n-3} u_x^5dx \nn\\
	& & {}-2(\gamma-\beta)(n-1) \io (x+\eps)^{-\beta+\gamma-1} g_\eps(x) u^{n-2} u_x^3 u_{xx}dx \nn\\
	& & {}-2(n-1) \io (x+\eps)^{-\beta+\gamma} g_{\eps x}(x) u^{n-2} u_x^3 u_{xx}dx \nn\\[2mm]
	&=:& I_{31}+ \cdots + I_{36}
	\qquad \mbox{for all } t\in (0,T).
  \eea
  In order to specify our choice of $\eps_\star$, 
  let us note that, according to our restriction on $n$ and with $P$ as specified in the formulation of
  the lemma, we have $P(n)>0$, which implies that when $n<2$,
  \bas
	& & \hspace*{-10mm}
	4(3-n) \Big\{ 6(n-1)-\frac{(4-n)^2}{4} - \frac{10(n-1)(2-n)}{3-n} \Big\} \\
	&=& 24(3-n)(n-1) - (3-n)(4-n)^2 - 40(n-1)(2-n) \\
	&=& -24 n^2 + 96n -72 -3n^2 + 24n - 48 + n^3 - 8n^2 +16n + 40n^2 -120n + 80 \\
	&=& P(n) >0.
  \eas
  Since in the case $n \in [2,3)$ we clearly have
  \bas
	6(n-1)-\frac{(4-n)^2}{4} \ge 6 (2-1) - \frac{(4-2)^2}{4}=5>0,
  \eas
  this entails that for any choice of $n\in (n_\star,3)$,
  \bas
	6(n-1) - \frac{(4-n)^2}{4} - \frac{10(n-1)(2-n)_+}{3-n} >0.
  \eas
  Consequently, with $\Lambda(\eps)$ as in Lemma \ref{lem1}, we can pick $\eps_\star\in (0,\eps_0)$ such that with
  $\Lambda_\star:=\Lambda(\eps_\star)$, we have
  \bas
	\Big\{ 6(n-1)-\frac{(4-n)^2}{4} \Big\} \Lambda_\star - \frac{10(n-1)(2-n)_+}{3-n}>0,
  \eas
  and thereupon fix a number $\mu\in (0,1)$ sufficiently close to $1$ and $\eta>0$ suitably small such that still
  \be{7.4}
	\Big\{6(n-1) - \frac{(4-n)^2}{4\mu} \Big\} \Lambda_\star - \frac{10(n-1)(2-n)_+}{3-n}
	- \Big\{ \Lambda_\star + \frac{50}{(3-n)^2} +1 \Big\} \eta >0,
  \ee
  and such that moreover
  \be{7.44}
	(1-\mu-\eta)\Lambda_\star - \eta>0.
  \ee
  Upon these choices, we first use Young's inequality to estimate $I_2$ according to
  \be{7.5}
	I_2 \le \mu \io (x+\eps)^{-\beta+\gamma} g_\eps(x) u^n u_{xxx}^2dx
	+ \frac{(4-n)^2}{4\mu} \io (x+\eps)^{-\beta+\gamma} g_\eps(x) u^{n-2} u_x^2 u_{xx}^2dx.
  \ee
  Next, recalling Lemma \ref{lem1}, we obtain
  \bea{7.6}
	I_4 &\le& \frac{\eta}{2} \io (x+\eps)^{-\beta+\gamma} g_\eps(x) u^n u_{xxx}^2dx
	+ c_1 \io (x+\eps)^{-\beta+\gamma} \frac{g_{\eps x}^2(x)}{g_\eps(x)} \cdot u^n u_{xx}^2dx \nn\\
	&\le& \frac{\eta}{2} |I_1| + c_2 \Gamma_1,
  \eea
  where
  \be{7.7}
	\Gamma_1 := \io (x+\eps)^{\alpha-\beta+\gamma-2} u^n u_{xx}^2dx
  \ee
  and $c_1$ and $c_2$, as all numbers $c_3,c_4,\ldots$ appearing below, denote positive constants depending on $n$, $\alpha$,
  $\beta$, and $\gamma$, but neither on $\eps \in (0,\eps_\star)$ nor on the solution $u$.\\
  Similarly, we find $c_3>0$ and $c_4>0$ such that
  \bea{7.8}
	I_5 &\le& \frac{\eta}{4} \io (x+\eps)^{-\beta+\gamma} g_\eps(x) u^n u_{xxx}^2 dx
	+ c_3 \io (x+\eps)^{-\beta+\gamma} \frac{g_{\eps x}^2(x)}{g_\eps(x)} \cdot u^{n-2} u_x^4dx \nn\\
	&\le& \frac{\eta}{4} |I_1| + c_4 \Gamma_2
  \eea
  with
  \be{7.9}
	\Gamma_2:=\io (x+\eps)^{\alpha-\beta+\gamma-2} u^{n-2} u_x^4dx,
  \ee
  and then $c_5>0$ and $c_6>0$ satisfying
  \bea{7.10}
	I_6 &\le& \frac{\eta}{8} \io (x+\eps)^{-\beta+\gamma} g_\eps(x) u^n u_{xxx}^2dx
	+ c_5 \io (x+\eps)^{-\beta+\gamma-2} g_\eps(x) u^n u_{xx}^2 dx \nn\\
	&\le& \frac{\eta}{8} |I_1| + c_5 \Gamma_1
  \eea
  and
  \bea{7.11}
	I_7 &\le& \frac{\eta}{2} \io (x+\eps)^{-\beta+\gamma} g_\eps(x) u^{n-2} u_x^2 u_{xx}^2dx
	+ c_6 \io (x+\eps)^{-\beta+\gamma-2} g_\eps(x) u^n u_{xx}^2dx \nn\\
	&\le& \frac{\eta}{2} \tilde I_{31} + c_6 \Gamma_1
  \eea
  where
  \be{7.12}
	\tilde I_{31}:=\frac{I_{31}}{-6(n-1)} = \io (x+\eps)^{-\beta+\gamma} g_\eps(x) u^{n-2} u_x^2 u_{xx}^2dx.
  \ee
  In much the same manner, we derive the inequalities
  \bea{7.13}
	I_8 &\le& \frac{\eta}{4} \io (x+\eps)^{-\beta+\gamma} g_\eps(x) u^{n-2} u_x^2 u_{xx}^2 dx
	+ c_7 \io (x+\eps)^{-\beta+\gamma-2} g_\eps(x) u^{n-2} u_x^4dx \nn\\
	&\le& \frac{\eta}{4} \tilde I_{31} + c_7 \Gamma_2
  \eea
  and
  \bea{7.14}
	I_{10} &\le& \frac{\eta}{8} \io (x+\eps)^{-\beta+\gamma} g_\eps(x) u^{n-2} u_x^2 u_{xx}^2 dx
	+ c_8 \io (x+\eps)^{-\beta+\gamma-2} \frac{g_{\eps x}^2(x)}{g_\eps(x)} \cdot u^n u_x^2dx \nn\\
	&\le& \frac{\eta}{8} \tilde I_{31} + c_9 \Gamma_3
  \eea
  for some positive $c_7, c_8$ and $c_9$ and
  \be{7.15}
	\Gamma_3:=\io (x+\eps)^{\alpha-\beta+\gamma-4} u^n u_x^2dx,
  \ee
  as well as
  \bea{7.16}
	I_{11} &\le& \frac{\eta}{16} \io (x+\eps)^{-\beta+\gamma} g_\eps(x) u^n u_{xxx}^2dx
	+ c_{10} \io (x+\eps)^{-\beta+\gamma-4} g_\eps(x) u^n u_x^2dx \nn\\
	&\le& \frac{\eta}{16} |I_1| + c_{10} \Gamma_3
  \eea
  and
  \bea{7.17}
	I_{12} &\le& \frac{\eta}{16} \io (x+\eps)^{-\beta+\gamma} g_\eps(x) u^{n-2} u_x^2 u_{xx}^2 dx
	+ c_{11} \io (x+\eps)^{-\beta+\gamma-4} g_\eps(x) u^n u_x^2dx \nn\\
	&\le& \frac{\eta}{16} \tilde I_{31} + c_{11} \Gamma_3
  \eea
  and
  \bea{7.18}
	I_{14} &\le& \frac{\eta}{32} \io (x+\eps)^{-\beta+\gamma} g_\eps(x) u^{n-2} u_x^2 u_{xx}^2dx
	+ c_{12} \io (x+\eps)^{-\beta+\gamma-4} \frac{g_{\eps x}^2(x)}{g_\eps(x)} \cdot u^{n+2}dx \nn\\
	&\le& \frac{\eta}{32} \tilde I_{31} + c_{13} \Gamma_4,
  \eea
  where 
  \be{7.19}
	\Gamma_4:=\io (x+\eps)^{\alpha-\beta+\gamma-6} u^{n+2}dx
  \ee
  and $c_{10}$, $c_{11}$, $c_{12}$, and $c_{13}$ are positive constants.\\
  As for the remaining terms on the right of (\ref{7.2}), we again apply Lemma \ref{lem1} to find $c_{14}>0$ and $c_{15}>0$
  such that
  \be{7.20}
	I_9 \le c_{14} \Gamma_1
  \ee
  and
  \be{7.21}
	I_{13} \le c_{15} J_2,
  \ee
  whereas Young's inequality provides $c_{16}>0$ fulfilling
  \bea{7.22}
	I_{15} &\le& c_{16} \io (x+\eps)^{\alpha-\beta+\gamma-3} u^{n-1} |u_x|^3dx \nn\\
	&\le& c_{16} \io (x+\eps)^{\alpha-\beta+\gamma-2} u^{n-2} u_x^4 dx
	+ c_{16} \io (x+\eps)^{\alpha-\beta+\gamma-6} u^{n+2}dx \nn\\
	&=& c_{16} \Gamma_2 + c_{16} \Gamma_4.
  \eea
  Finally, abbreviating
  \be{7.23}
	\tilde I_{32}:=\io (x+\eps)^{-\beta+\gamma} g_\eps(x) u^{n-4} u_x^6dx,
  \ee
  using Young's inequality we obtain constants $c_{17},...,c_{21}$ such that
  \bea{7.24}
	I_{33} &\le& \frac{\eta}{2} \io (x+\eps)^{-\beta+\gamma} g_\eps(x) u^{n-4} u_x^6dx
	+ c_{17} \io(x+\eps)^{-\beta+\gamma-2} g_\eps(x) u^{n-2} u_x^4dx \nn\\
	&\le& \frac{\eta}{2} \tilde I_{32} + c_{17} \Gamma_2
  \eea
  and
  \bea{7.25}
	I_{34} &\le& \frac{\eta}{4} \io (x+\eps)^{-\beta+\gamma} g_\eps(x) u^{n-4} u_x^6dx
	+ c_{18} \io (x+\eps)^{-\beta+\gamma} \frac{g_{\eps x}^2(x)}{g_\eps(x)} \cdot u^{n-2} u_x^4dx \nn\\
	&\le& \frac{\eta}{4} \tilde I_{32} + c_{18} \Gamma_2
  \eea
  as well as
  \bea{7.26}
	I_{35} &\le& \frac{\eta}{8} \io (x+\eps)^{-\beta+\gamma} g_\eps(x) u^{n-4} u_x^6dx
	+ c_{19} \io (x+\eps)^{-\beta+\gamma-2} g_\eps(x) u^n u_{xx}^2dx \nn\\
	&\le& \frac{\eta}{8} \tilde I_{32} + c_{19} \Gamma_1
  \eea
  and
  \bea{7.27}
	I_{36} &\le& \frac{\eta}{16} \io (x+\eps)^{-\beta+\gamma} g_\eps(x) u^{n-4} u_x^6dx
	+ c_{20} \io (x+\eps)^{-\beta+\gamma} \frac{g_{\eps x}^2(x)}{g_\eps(x)} \cdot u^n u_{xx}^2 dx \nn\\
	&\le& \frac{\eta}{16} \tilde I_{32} + c_{21} \Gamma_1.
  \eea
  In light of (\ref{7.5})-(\ref{7.27}), (\ref{7.2}) and (\ref{7.3}) thus yield
  \bea{7.28}
	\frac{1}{2} \frac{d}{dt} \io (x+\eps)^\gamma u_x^2dx
	&\le& - (1-\mu-\eta) \io (x+\eps)^{-\beta+\gamma} g_\eps(x) u^n u_{xxx}^2dx \nn\\
	& & {}- \Big\{ 6(n-1)-\frac{(4-n)^4}{4\mu} - \eta \Big\} \cdot \io (x+\eps)^{-\beta+\gamma} g_\eps(x) u^{n-2} u_x^2
		u_{xx}^2dx \nn\\
	& & {}+ \Big\{ - \frac{2}{5} (n-1)(n-2)(3-n) + \eta \Big\} \cdot \io (x+\eps)^{-\beta+\gamma} g_\eps(x) u^{n-4} u_x^6dx
		\nn\\
	& & {}+ c_{22} \io (x+\eps)^{\alpha-\beta+\gamma-2} u^n u_{xx}^2 dx
	+ c_{22} \io (x+\eps)^{\alpha-\beta+\gamma-2} u^{n-2} u_x^4dx \nn\\
	& & {}+ c_{22} \io (x+\eps)^{\alpha-\beta+\gamma-4} u^n u_x^2dx
	+ c_{22} \io (x+\eps)^{\alpha-\beta+\gamma-6} u^{n+2}dx
  \eea
  for all $t\in (0,T)$ with some $c_{22}>0$, where we have used that 
  $\sum_{j=1}^N \frac{\eta}{2^j} < \eta$ for all $N\in\N$.\abs
  Now by means of Lemma \ref{lem6} we can find $c_{23}>0$ such that
  \bas
	& & \hspace*{-20mm}
	\Big\{-\frac{2}{5}(n-1)(n-2)(3-n)+\eta \Big\} \io (x+\eps)^{-\beta+\gamma} g_\eps(x) u^{n-4} u_x^6dx \\
	&\le& \Big\{\frac{2}{5}(n-1)(2-n)_+(3-n)+\eta\Big\} \io (x+\eps)^{\alpha-\beta+\gamma} u^{n-4} u_x^6dx \nn\\
	&\le& \frac{25}{(1-\eta)(3-n)^2} 
	\Big\{ \frac{2}{5}(n-1)(2-n)_+(3-n)+\eta \Big\}
	\io (x+\eps)^{\alpha-\beta+\gamma} u^{n-2} u_x^2 u_{xx}^2dx \nn\\
	& & {}+ c_{23} \io (x+\eps)^{\alpha-\beta+\gamma-2} u^{n-2} u_x^4dx \\
	&\le& \frac{25}{(3-n)^2} \Big\{ \frac{2}{5}(n-1)(2-n)_+(3-n)+2\eta \Big\}
		\io (x+\eps)^{\alpha-\beta+\gamma} u^{n-2} u_x^2 u_{xx}^2dx \nn\\
	& & {}+ c_{23} \io (x+\eps)^{\alpha-\beta+\gamma-2} u^{n-2} u_x^4dx.
  \eas
	The last inequality follows from the fact that for
	$A:=\frac25(n-1)(2-n)_+(3-n)<1$ ($0<n<3$), we have $(A+\eta)/(1-\eta)\le
	A+2\eta$ if $0<\eta<\frac12(1-A)$.
  Then applying Lemma \ref{lem_inter}, we obtain $c_{24}>0$ satisfying
  \bea{7.66}
	c_{22} \io (x+\eps)^{\alpha-\beta+\gamma-2} u^n u_{xx}^2dx
	&+& (c_{22}+c_{23}) \io (x+\eps)^{\alpha-\beta+\gamma-2} u^{n-2} u_x^4 dx \nn \\
	&& {}+ c_{22} \io (x+\eps)^{\alpha-\beta+\gamma-4} u^n u_x^2dx \nn\\
	&\le& \eta \io (x+\eps)^{\alpha-\beta+\gamma} u^n u_{xxx}^2 dx
	+ \eta \io (x+\eps)^{\alpha-\beta+\gamma} u^{n-2} u_x^2 u_{xx}^2dx \nn\\
	& & {}+ c_{24} \io (x+\eps)^{\alpha-\beta+\gamma-6} u^{n+2}dx.
  \eea
  Therefore, (\ref{7.28}) shows that
  \bas
	&&\frac{1}{2} \frac{d}{dt} \io (x+\eps)^\gamma u_x^2dx
	\le -(1-\mu-\eta) \io (x+\eps)^{-\beta+\gamma} g_\eps(x) u^n u_{xxx}^2dx \\
	& & {}- \Big\{ 6(n-1) - \frac{(4-n)^2}{4\mu} - \eta \Big\} \cdot \io (x+\eps)^{-\beta+\gamma} g_\eps(x)
		u^{n-2} u_x^2 u_{xx}^2dx \\
	& & {}+ \eta \io (x+\eps)^{\alpha-\beta+\gamma} u^n u_{xxx}^2dx \nn\\
	& & {}+ \bigg\{ \frac{25}{(3-n)^2} \Big[ \frac{2}{5}(n-1)(2-n)_+(3-n) +2\eta\Big] + \eta \bigg\}
		\io (x+\eps)^{\alpha-\beta+\gamma} u^{n-2} u_x^2 u_{xx}^2dx \\
	& & {}+ (c_{22}+c_{24}) \io (x+\eps)^{\alpha-\beta+\gamma-6} u^{n+2}dx
	\qquad \mbox{for all } t\in (0,T).
  \eas
  Since clearly $1-\mu-\eta$ and $6(n-1)-\frac{(4-n)^2}{4\mu}-\eta$ are both positive thanks to (\ref{7.44}) and
  (\ref{7.4}), we may now use the lower estimate for $g_\eps$ established in Lemma \ref{lem1} to infer that
  \bas
	\frac{1}{2} \frac{d}{dt} \io (x+\eps)^\gamma u_x^2dx
	&\le& - \Big\{ (1-\mu-\eta)\Lambda_\star - \eta \Big\} \io (x+\eps)^{\alpha-\beta+\gamma} u^n u_{xxx}^2dx \\
	& & {}- \bigg\{ \Big[6(n-1)-\frac{(4-n)^2}{4\mu}-\eta\Big] \Lambda_\star
		-\frac{10(n-1)(2-n)_+}{3-n} - \frac{50\eta}{(3-n)^2} - \eta \bigg\} \\
	& & \hspace*{7mm}\times \io (x+\eps)^{\alpha-\beta+\gamma} u^{n-2} u_x^2 u_{xx}^2 dx \\
	& & {}+ (c_{22}+c_{24}) \io (x+\eps)^{\alpha-\beta+\gamma-6} u^{n+2}dx
	\qquad \mbox{for all } t\in (0,T),
  \eas
  because $\eps<\eps_\star$ and hence $\Lambda(\eps)\ge \Lambda_\star$ by the monotonicity of $\Lambda$ asserted
  by Lemma \ref{lem1}. According to (\ref{7.44}) and
  (\ref{7.4}), after another application of (\ref{7.66}), this entails (\ref{7.1}).
\qed
Under additional assumptions on the parameters $\alpha$, $\beta$, and $\gamma$,
we are able to derive a priori estimates for small times
only depending on the initial data.
More precisely, if the parameter $\gamma$ is chosen large enough, then the weight in the integral on the right-hand side
of (\ref{7.1}) is sufficiently regular, whence from the above we can deduce a bound for
$\io (x+\eps)^\gamma u_x^2(x,t)dx$ for all sufficiently small $t>0$. 
Since we plan to finally achieve a boundedness property for $u$ itself with respect to the norm in 
$L^\infty(\Omega)$, we require that $\gamma<1$. 
This explains the restriction on $\beta$ in the following lemma.
\begin{lem}[A priori estimate for small times]\label{lem13}
  Let $n_*=1.5361\ldots$ and $\eps_\star\in (0,1)$ be as in Lemma \ref{lem7},
  let $\alpha>0$, $\beta\in (-1,\alpha-4)$, and 
  \be{13.1}
	\gamma\in(5-\alpha+\beta,1).
  \ee
  Then one can find $c>0$ such that
  for all $A>0$ and $B>0$, there exists $T_0(A,B) \in (0,1)$ with the following property:
  If for some $\eps \in (0,\eps_\star)$ and $T\in (0,T_0(A,B))$, $u\in C^{4,1}(\bar\Omega\times [0,T))$
  is positive and solves (\ref{0eps}) in $\Omega \times (0,T)$ with
  \be{13.23}
	\io (x+\eps)^\gamma u_x^2 (x,0)dx \le A
	\qquad \mbox{and} \qquad
	\io (x+\eps)^\beta u(x,0) dx \le B,
  \ee
  then
  \bea{13.4}
	&&\hspace*{-10mm}\sup_{t\in (0,T)} \io (x+\eps)^\gamma u(x,t)dx 
	\le c\int_0^T \io (x+\eps)^{\alpha-\beta+\gamma} u^n u_{xxx}^2dxdt \\
	&+& c\int_0^T \io (x+\eps)^{\alpha-\beta+\gamma} u^{n-2} u_x^2 u_{xx}^2dxdt
	+ c\int_0^T \io (x+\eps)^{\alpha-\beta+\gamma} u^{n-4} u_x^6dxdt \nn\\
	&+& c\int_0^T \io (x+\eps)^{\alpha-\beta+\gamma-2} u^n u_{xx}^2 dxdt
	+ c\int_0^T \io (x+\eps)^{\alpha-\beta+\gamma-2} u^{n-2} u_x^4dxdt \nn\\[2mm]
	&\le& A+1.
  \eea
  In particular, in that case 
  there exists $C(A,B)>0$ such that the flux $J$, defined in (\ref{J}), satisfies
  \be{13.44}
	\int_0^T \io (x+\eps)^{-\alpha-\beta+\gamma} J_x^2 dxdt \le C(A,B).
  \ee
\end{lem}
\proof
  Let us first note that our hypothesis $\beta\in (-1,\alpha-4)$ entails the inequality
  $5-\alpha+\beta<1$, whence the assumption $\gamma\in (5-\alpha+\beta,1)$ indeed is meaningful. 
  Then with $T_0(A,B) \in (0,1)$ to be fixed below, we assume that $T\in (0,T_0(A,B))$ and that $u$ has the 
  properties listed above. 
  Thus, for each $t\in (0,T)$, by (\ref{13.23}) and Lemma \ref{lem8}, we have
  $\io (x+\eps)^\beta u(x,t) dx \le B$, so that Lemma \ref{lem32} says that
  \be{13.444}
	u(x,t) \le c_1 B + c_1 \Big( \io (x+\eps)^\gamma u_x^2(x,t)dx \Big)^\frac{1}{2}
	\qquad \mbox{for all } x\in\Omega
  \ee
  with some $c_1>0$, where we have used that $\gamma<1$.
  Consequently, thanks to (\ref{13.1}), the integral on the right-hand side of (\ref{7.1}) 
  can be estimated according to
  \bas
	\io (x+\eps)^{\alpha-\beta+\gamma-6} u^{n+2}(x,t)dx
	&\le& 2^{n+2} (c_1 B)^{n+2} \Big( \io (x+\eps)^{\alpha-\beta+\gamma-6} dx \Big)^{n+2} \nn\\
	& & {}+ 2^{n+2} \cdot L \cdot c_1^{n+2} \Big( \io (x+\eps)^\gamma u_x^2(x,t)dx \Big)^\frac{n+2}{2} \nn\\
	&\le& c_2(B) + c_3 \Big( \io (x+\eps)^\gamma u_x^2(x,t)dx \Big)^\frac{n+2}{2}
  \eas
  with appropriate constants $c_2(B)>0$ and $c_3>0$. 
  From Lemmas \ref{lem7} and \ref{lem6},
	we thus obtain $c_4>0$, $c_5(B)>0$, and $c_6>0$ such that
  \bea{13.5}
	\frac{d}{dt} \io (x+\eps)^\gamma u_x^2dx
	&+& c_4 \io (x+\eps)^{\alpha-\beta+\gamma} u^n u_{xxx}^2dx
	+ c_4 \io (x+\eps)^{\alpha-\beta+\gamma} u^{n-2} u_x^2 u_{xx}^2dx \nn\\
	&+& c_4 \io (x+\eps)^{\alpha-\beta+\gamma} u^{n-4} u_x^6dx \nn\\
	&+& c_4  \io (x+\eps)^{\alpha-\beta+\gamma-2} u^n u_{xx}^2 dx
	+ c_4 \io (x+\eps)^{\alpha-\beta+\gamma-2} u^{n-2} u_x^4dx \nn\\
	&\le& c_5(B) + c_6 \Big( \io (x+\eps)^\gamma u_x^2 \Big)^\frac{n+2}{2}dx
	\qquad \mbox{for all } t\in (0,T).
  \eea
  With the above constants being fixed, we consider the solution $y\equiv y_{A,B}$ of the initial-value problem
  \bas
	\left\{ \begin{array}{l}
	y'(t) = c_5(B) + c_6 y^\frac{n+2}{2}(t), \qquad t>0, \\[1mm]
	y(0)=A.
	\end{array} \right.
  \eas
  It is then clearly possible to fix some sufficiently small $T_0(A,B)\in (0,1)$
  such that $y(t) \le A+1$ for all $t\in (0,T_0(A,B))$,
  and a comparison argument for ordinary differential equations,
	applied to (\ref{13.5}), shows that
  \bas
	\io (x+\eps)^\gamma u_x^2(x,t) dx \le A+1
	\qquad \mbox{for all } t\in (0,T), \ T<T_0(A,B).
  \eas
  Inserting this into (\ref{13.5}) and integrating, we readily arrive at (\ref{13.4}).\\
  From this, the estimate (\ref{13.44}) easily follows upon recalling (\ref{Jx}), (\ref{13.444}), and 
  Lemma \ref{lem1} and applying Lemma \ref{lem32} and Lemma \ref{lem_inter}.
\qed
\mysection{Local existence in the approximate problems}\label{sec.locex}
The a priori estimate of Lemma \ref{lem13} allows us to prove the
local existence of a classical solution to the approximate problem (\ref{0eps})
for smooth initial data $u_0$ with compactly supported derivative $u_{0x}$.
\begin{lem}[Local existence for smooth data]\label{lem9}
  Let $\eps_0=\min\{1,\sqrt{\frac{L}{2}}\}$ and $\eps\in(0,\eps_0)$, and let 
  $u_0\in C^\infty(\overline\Omega)$ be positive and such that 
  $u_{0x}\in C_0^\infty(\Omega)$. Then there exist $\tm\in(0,\infty]$
  and a unique positive classical solution $u\in C^{4,1}(\bar\Omega \times [0,\tm))$
  of (\ref{0eps}) in $\Omega \times (0,\tm)$
  with $u(x,0)=u_0(x)$ for all $x\in \Omega$. Moreover, $\tm$ has the property that
  \be{9.1}
	\mbox{if $\tm<\infty$ then either } \liminf_{t\nearrow \tm} \Big(\inf_{x\in\Omega} u(x,t)\Big)=0 
	\quad \mbox{ or } \limsup_{t\nearrow \tm} \Big(\sup_{x\in\Omega} u(x,t)\Big) =\infty.
  \ee
\end{lem}
\proof
  For $k\in\N$, we let $f_k\in C^\infty(\R)$ be a smooth nondecreasing truncation function on $\R$ such 
  that $f_k(s)=s$ for all $s\in [\frac{1}{k},k]$
  and $\frac{1}{2k} \le f_k \le 2k$ on $\R$. Then each of the problems
  \be{0epsk}
	\left\{ \begin{array}{l}
	u_{kt} = \frac{1}{(x+\eps)^\beta} \cdot \Big\{ -g_\eps(x) f_k^n(u_k) u_{kxx} 
	+ 2g_\eps(x) f_k^{n-1}(u_k) u_{kx}^2 \Big\}_{xx}, \qquad x\in \Omega, \ t>0, \\[1mm]
	u_{kx}=u_{kxxx}=0, \qquad x\in\pO, \ t>0, \\[1mm]
	u_k(x,0)=u_0(x), \qquad x\in\Omega,
	\end{array} \right.
  \ee
  is non-degenerate, and since $u_{0x}$ has compact support in $\Omega$, standard parabolic theory \cite{Fri64}
  yields a uniquely determined global solution $u\in C^{4.1}(\bar\Omega \times [0,\infty))$.\\
  Now for sufficiently large $k_0\in\N$ and each $k\ge k_0$, it follows from the continuity of $u_k$ and the
  positivity of $u_0$ in $\bar\Omega$ that
  \bas
	T_k:=\sup \Big\{ T>0 \ \Big| \ \frac{1}{k} \le u_k \le k \mbox{ in } \Omega \times (0,T) \Big\}
  \eas
  is a well-defined element of $(0,\infty]$, and by uniqueness in (\ref{0epsk}), it is clear that the sequence
  $(T_k)_{k\ge k_0}$ is nondecreasing, and that $u_{k_2}\equiv u_{k_1}$ in $\Omega \times (0,T_{k_1})$ whenever
  $k_2 \ge k_1 \ge k_0$.
  Consequently, the definition $\tm:=\lim_{k\to\infty} T_k \in (0,\infty]$ is meaningful, 
  and the trivially existing pointwise limit $u(x,t):=\lim_{k\to\infty} u_k(x,t)$, $(x,t)\in \bar\Omega \times [0,\tm)$,
  satisfies $u\equiv u_k$ in $\Omega \times (0,T_k)$ for each $k\ge k_0$. 
  It is therefore evident from (\ref{0epsk}) and the definition of $f_k$ 
  that $u$ actually solves (\ref{0eps}) in $\Omega \times (0,\tm)$ with $u|_{t=0}=u_0$ in $\Omega$.\\
  It remains to verify (\ref{9.1}).
  To this end, we assume on the contrary that $\tm<\infty$, but that both
  $\liminf_{t\nearrow \tm} (\inf_{x\in\Omega} u(x,t))>0$ and 
  $\limsup_{t\nearrow\tm} (\sup_{x\in\Omega} u(x,t))<\infty$. 
  Then for some $k\ge k_0$, we would have $\frac{2}{k} \le u \le \frac{k}{2}$ in $\Omega \times (0,\tm)$, implying that
  $u\equiv u_k$ in $\Omega \times (0,\tm)$ by uniqueness. 
  But since $u_k$ is continuous at $t=\tm$, this would entail that $T_k>\tm$ and hence contradict the definition
  of $\tm$.
\qed
The following result rules out the occurrence of the first alternative in (\ref{9.1}); that is,
solutions to the approximate problem (\ref{0eps}) cannot develop a dead core within finite time.
\begin{lem}[Absence of dead core formation]\label{lem12}
  Let $n>1$, $\alpha> 0$, and $\beta\in\R$. 
  Then for all $\eps \in (0,\eps_0), \delta>0, M>0$, and $T>0$ there exists
  $C(\eps,\delta,M,T)>0$ such that if $u\in C^{4,1}(\bar\Omega\times [0,T))$ is a positive classical solution
  of (\ref{0eps}) in $\Omega\times (0,T)$ satisfying
  \be{12.01}
	u(x,0) \ge \delta
	\qquad \mbox{for all } x\in\Omega
  \ee
  and
  \be{12.1}
	u(x,t) \le M
	\qquad \mbox{for all $x\in\Omega$ and } t\in (0,T),
  \ee
  we have the inequality
  \be{12.2}
	\io \frac{1}{u^2(x,t)}dx \le C(\eps,\delta,M,T)
	\qquad \mbox{for all } t\in (0,T).
  \ee
\end{lem}
\proof
  Our goal is to conclude (\ref{12.2}) from a differential inequality for $\int_\Omega(x+\eps)^\beta u^{-2}$
  which we shall thus derive first. 
  To this end, we twice integrate by parts over $\Omega$ to compute, using $J$ as defined in (\ref{J}),
  \bas
	\frac{d}{dt} \io (x+\eps)^\beta  \frac{1}{u^2}dx
	&=& - 2\io (x+\eps)^\beta \frac{u_t}{u^2}dx
	= - 2\io \frac{1}{u^3} J_{xx} dx \\
	&=& - 6 \io \frac{u_x}{u^4} J_x dx
	= 6 \io \Big(\frac{u_x}{u^4}\Big)_x Jdx \nn\\
	&=& 6 \io \Big[ \frac{u_{xx}}{u^4} - 4\frac{u_x^2}{u^5} \Big]
	\Big[ -g_\eps(x) u^n u_{xx} + 2g_\eps(x) u^{n-1} u_x^2 \Big]dx \nn\\
	&=& - 6\io g_\eps(x) u^{n-4} u_{xx}^2dx
	+ 36 \io g_\eps(x) u^{n-5} u_x^2 u_{xx}dx
	- 48 \io g_\eps(x) u^{n-6} u_x^4dx
  \eas
  for all $t\in (0,T)$, because $u_x=J_x=0$ on $\pO$ according to Lemma \ref{lem_J}.
  Since one more integration by parts yields
  \bas
	36\io g_\eps(x) u^{n-5} u_x^2 u_{xx}dx
	&=& 12(5-n) \io g_\eps(x) u^{n-6} u_x^4dx
	- 12 \io g_{\eps x}(x) u^{n-5} u_x^3dx,
  \eas
  this shows that
  \bea{12.4}
	\frac{d}{dt} \io (x+\eps)^\beta \cdot \frac{1}{u^2}dx
	&=& - 6\io g_\eps(x) u^{n-4} u_{xx}^2dx
	-12(n-1) \io g_\eps(x) u^{n-6} u_x^4dx \nn \\
	&&{}-12 \io g_{\eps x}(x) u^{n-5} u_x^3dx
  \eea
  for all $t\in (0,T)$.
  Here, since $n>1$, the second term on the right-hand side is nonpositive, and
  by means of Young's inequality and Lemma \ref{lem1}, we can find $c_1>0$ and $c_2>0$ fulfilling
  \bas
	-12 \io g_{\eps x}(x) u^{n-5} u_x^3dx
	&\le& 12(n-1) \io g_\eps(x) u^{n-6} u_x^4dx
	+ c_1 \io \frac{g_{\eps x}^4(x)}{g_\eps^3(x)} \cdot u^{n-2} dx\\
	&\le& c_2 \io (x+\eps)^{\alpha-4} u^{n-2}dx,
  \eas
  whence (\ref{12.4}) in particular entails that
  \be{12.5}
	\frac{d}{dt} \io (x+\eps)^\beta \cdot \frac{1}{u^2}dx
	\le c_2 \io (x+\eps)^{\alpha-4} u^{n-2}dx
	\qquad \mbox{for all } t\in (0,T).
  \ee
  Now if $n\ge 2$, writing $c_3(\eps):=\io (x+\eps)^{\alpha-4}dx$ and using (\ref{12.1}), from (\ref{12.5}) we obtain
  \bas
	\frac{d}{dt} \io (x+\eps)^\beta \cdot \frac{1}{u^2}dx
	\le c_2 c_3(\eps) M^{n-2}
	\qquad \mbox{for all } t\in (0,T),
  \eas
  which after integration implies that
  \bas
	\io (x+\eps)^\beta \cdot \frac{1}{u^2(x,t)} dx
	\le \io (x+\eps)^\beta \cdot \frac{1}{u^2(x,0)} dx
	+ c_2 c_3(\eps) M^{n-2} T
	\qquad \mbox{for all } t\in (0,T).
  \eas
  As a consequence of (\ref{12.01}), we thereby find that
  \bas
	\io \frac{1}{u^2(x,t)} dx
	&\le& c_4(\eps) \io (x+\eps)^\beta \cdot \frac{1}{u^2(x,t)} dx \\
	&\le& c_4(\eps) \Big\{ \frac{c_5(\eps)}{\delta^2} + c_2 c_3(\eps) M^{n-2}T \Big\}
	\qquad \mbox{for all } t\in (0,T)
  \eas
  with $c_4(\eps):=\max\{\eps^{-\beta}, (L+\eps)^{-\beta}\}$ and $c_5(\eps):=\io (x+\eps)^\beta dx$.\\
  In the remaining case $n<2$, we first apply the H\"older inequality 
	with $p=\frac{2}{n}$ and $p'=\frac{2}{2-n}$ to the 
	right-hand side in (\ref{12.5}) to obtain
  \bas
	\frac{d}{dt} \io (x+\eps)^\beta \cdot \frac{1}{u^2}dx
	&\le& c_2\io (x+\eps)^{\alpha-4-(2-n)\frac{\beta}{2}}
	((x+\eps)^\beta u^{-2})^{\frac{2-n}{2}}dx \\
	&\le& c_2 c_6(\eps) \Big( \io (x+\eps)^\beta \cdot \frac{1}{u^2}dx 
	\Big)^\frac{2-n}{2}
	\qquad \mbox{for all } t\in (0,T)
  \eas
  with
  \bas
	c_6(\eps) := \Big( \io (x+\eps)^\frac{2(\alpha-4)-(2-n)\beta}{n} dx \Big)^\frac{n}{2}.
  \eas
  Integrating this in time shows that in this case,
  \bas
	\io (x+\eps)^\beta \cdot \frac{1}{u^2(x,t)} dx
	\le \Big\{ \io (x+\eps)^\beta \cdot \frac{1}{u^2(x,0)} dx
	+ \frac{n}{2} c_2 c_6(\eps)T \Big\}^\frac{2}{n}
	\qquad \mbox{for all } t\in (0,T),
  \eas
  and hence,
  \bas
	\io \frac{1}{u^2(x,t)} dx
	\le c_4(\eps) \cdot \Big\{ \frac{c_5(\eps)}{\delta^2}
	+ \frac{n}{2} c_2 c_6(\eps) T \Big\}^\frac{2}{n}
	\qquad \mbox{for all } t\in (0,T),
  \eas
  according to (\ref{12.1}).
\qed
\mysection{H\"older continuity}\label{sec.hoelder}
We next derive a spatio-temporal H\"older estimate for the above solutions to the approximate problems.
This will allow us to construct a continuous weak solution of (\ref{0}) along a uniformly convergent sequence
of appropriate solutions of (\ref{0eps}) as $\eps\to 0$.
\begin{lem}[H\"older estimate]\label{lem14}
  With $n_\star$ as in Lemma \ref{lem7},
  assume that $n\in (n_\star,3)$ and that $\alpha>0$, $\beta<\alpha-4$, and 
  $\gamma<1$ are such that $\alpha-\beta+\gamma>5$. 
  Moreover, let $A>0$ and $B>0$, and let $\eps_\star$ and $T_0(A,B) \in (0,1)$ be as given by Lemma \ref{lem7} and 
  Lemma \ref{lem13}, respectively. 
  Then there exists $C(A,B)>0$ such that, whenever $u\in C^{4,1}(\bar\Omega \times [0,T))$ is a positive classical solution
  of (\ref{0eps}) in $\bar \Omega \times (0,T)$ for some $T\in (0,T_0(A,B))$ with
  \bas
	\io (x+\eps)^\gamma u_x^2(x,0)dx \le A
	\qquad \mbox{and} \qquad
	\io (x+\eps)^\beta u(x,0)dx \le B,
  \eas
  the estimate
  \be{14.1}
	|u(x_2,t_2)-u(x_1,t_1)| \le C(A,B) \cdot \Big(|x_2-x_1|^\theta + |t_2-t_1|^\frac{\theta}{2\theta+3} \Big)
  \ee
  holds for all $x_1$, $x_2\in \Omega$ and $t_1$, $t_2\in (0,T)$ with $\theta:=\min\{\frac{1}{2}, \frac{1-\gamma}{2}\}$.
\end{lem}
\proof
  According to Lemma \ref{lem13}, we can pick $c_1$, as well as all 
	constants $c_2,\ldots$ below, possibly depending on
  $n$, $\alpha$, $\beta$, $\gamma$, $A$, $B$, and $L$ but independent from 
	$\eps$ and $u$, such that 
  \bas
	\io (x+\eps)^\gamma u_x^2 \le c_1
	\qquad \mbox{for all  } t\in (0,T).
  \eas
  Hence, Lemma \ref{lem31} provides $c_2>0$ such that we have the spatial H\"older estimate
  \be{14.2}
	|u(x_2,t_0)-u(x_1,t_0)| \le c_2 |x_2-x_1|^\theta
	\qquad \mbox{for all $x_1,x_2\in\Omega$ and } t_0 \in (0,T).
  \ee
  Using this, a corresponding H\"older estimate with respect to the time variable, that is, the inequality
  \be{14.3}
	|u(x_0,t_2)-u(x_0,t_1)| \le M |t_2-t_1|^\frac{\theta}{2\theta+3}
	\qquad \mbox{for all $x_0\in\Omega$ and } t_1,t_2 \in (0,T)
  \ee
  with suitably large $M>1$, can be derived by adapting a standard technique due to Gilding and Kru\v{z}kov 
  (\cite{Gil76}, cf.~also \cite{BeFr90} for a related procedure in a fourth-order setting).
  Indeed, following \cite{BeFr90}, let us assume that (\ref{14.3}) be false, meaning that
  for some $x_0\in\Omega$ and $t_1,t_2\in (0,T)$ we have
  \be{14.33}
	u(x_0,t_2)-u(x_0,t_1)>M |t_2-t_1|^\frac{\theta}{2\theta+3},
  \ee
  where for definiteness we may suppose that $t_1<t_2$.
  We then fix any $\zeta \in C_0^\infty(\R)$ such that $0\le \zeta \le 1$ on $\R$, $\zeta \equiv 1$ in
  $[-\frac{1}{2},\frac{1}{2}]$ and $\zeta \equiv 0$ in $\R \setminus [-1,1]$, and let
  \bas
	\psi(x):=\zeta \Big(\frac{x-x_0}{\eta}\Big), \qquad x\in\bar\Omega,
  \eas
  with  
  \be{eta}
	\eta:=\Big(\frac{M}{16 c_2} \Big)^\frac{1}{\theta} (t_2-t_1)^\frac{1}{2\theta+3}.
  \ee
  Furthermore, we introduce the functions $\xi_\delta$, $\delta>0$, given by
  \be{14.4}
	\xi_\delta(t):=\frac{1}{\delta} \int_{-\infty}^t \Big\{ \zeta\Big(\frac{s-t_2}{\delta}\Big)
	- \zeta\Big(\frac{s-t_1}{\delta}\Big)\Big\} ds, \qquad t\in (0,T),
  \ee
  which belong to $C_0^\infty((0,T))$ and satisfy $0\ge \xi_\delta \ge -c_3$ with 
  $c_3:=\int_{-1}^t \zeta(\sigma) d\sigma$, provided that $\delta<\delta_0:=\min\{t_1,T-t_2\}$ (this ensures that $\xi_\delta(0)=\xi_\delta(T)=0$).
  Therefore, testing (\ref{0eps}) against $\psi(x)\xi_\delta(t)$, $(x,t)\in\Omega\times (0,T)$, we obtain
  \be{14.5}
	\int_0^T \io u(x,t) \psi(x)\xi_\delta'(t) dxdt
	= \int_0^T \io \Big[ (x+\eps)^{-\beta} \psi(x)\xi_\delta(t) \Big]_x J_x(x,t) dxdt
	\qquad \mbox{for all } \delta\in (0,\delta_0)
  \ee
  with $J$ as defined in (\ref{J}), where we again have used that $J_x=0$ on $\pO$
  by (\ref{Jx}).\\
  We insert the definition of $\xi'_\delta(t)$, substitute $\sigma=\frac{t-t_i}{\delta}$, $i\in \{1,2\}$,
  and perform the limit $\delta\to 0$, to estimate the left-hand side in (\ref{14.5}) from below according to
  \bas
	\frac{1}{c_3} \lim_{\delta\searrow 0} \int_0^T \io u(x,t)\psi(x)\xi_\delta'(t) dxdt
	&=& \frac{1}{c_3} \lim_{\delta\searrow 0} \io \int_{-1}^1 
		\Big[ u(x,t_2+\delta\sigma)-u(x,t_1+\delta\sigma) \Big] \cdot \zeta(\sigma) d\sigma
		\cdot \psi(x)dx \nn\\
	&=& \io \Big[u(x,t_2)-u(x,t_1)\Big] \cdot \psi(x) dx \nn\\
	&\ge& \io \Big[u(x_0,t_2)-u(x_0,t_1)\Big] \cdot \psi(x) dx \nn\\
	& & - \io \Big[ |u(x,t_2)-u(x_0,t_2)| + |u(x_0,t_1)-u(x,t_1)| \Big] \cdot \psi(x) dx,
  \eas
  whence using (\ref{14.33}) and (\ref{14.2}) yields
  \bea{14.6}
	\frac{1}{c_3} \lim_{\delta\searrow 0} \int_0^T \io u(x,t)\psi(x)\xi_\delta'(t) dxdt
	&\ge& M (t_2-t_1)^\frac{\theta}{2\theta+3} \io \psi(x) dx
	- 2c_2 \io |x-x_0|^\theta \cdot \psi(x)dx \nn\\
	&\ge& M(t_2-t_1)^\frac{\theta}{2\theta+3} \cdot \frac{\eta}{2}
	-2c_2 \cdot \eta^\theta \cdot 2\eta \nn\\
	&=& \frac{\eta}{2} \cdot \Big\{ M(t_2-t_1)^\frac{\theta}{2\theta+3} - 8c_2 \eta^\theta \Big\} \nn\\
	&=& \frac{\eta}{2} \cdot \frac{1}{2} M (t_2-t_1)^\frac{\theta}{2\theta+3} \nn\\
	&=& c_4 M^{1+\frac{1}{\theta}} (t_2-t_1)^\frac{\theta+1}{2\theta+3}
  \eea  
  with $c_4:=[4\cdot (16c_2)^\frac{1}{\theta}]^{-1}$.
  On the right-hand side of (\ref{14.5}), by the Cauchy-Schwarz inequality and Lemma \ref{lem13},
  we find $c_5>0$ fulfilling
  \bea{14.7}
	& & \hspace*{-10mm}
	\bigg| \int_0^T \io \Big[ (x+\eps)^{-\beta} \psi(x)\xi_\delta(t) \Big]_x \cdot J_x dxdt\bigg| \nn\\
	&\le& \bigg( \int_0^T \io (x+\eps)^{-\alpha-\beta+\gamma} J_x^2dxdt \bigg)^\frac{1}{2} 
	\bigg( \int_0^T \xi_\delta^2(t) dt \bigg)^\frac{1}{2} 
	\bigg( \io (x+\eps)^{\alpha+\beta-\gamma} \Big[(x+\eps)^{-\beta} \psi(x) \Big]_x^2 dx \bigg)^\frac{1}{2}
		\nn\\
	&\le& c_5 c_3 (t_2-t_1+2\delta)^\frac{1}{2} \cdot 
	\bigg( \io (x+\eps)^{\alpha+\beta-\gamma} \Big[(x+\eps)^{-\beta} \psi(x) \Big]_x^2 dx \bigg)^\frac{1}{2},
  \eea
  since $\xi_\delta(t)=0$ for $t\le t_1-\delta$ or $t\ge t_2+\delta$ and
	$\xi_\delta^2\le c_3^2$.
	We use $\alpha-\beta-\gamma>5-2\gamma>3$ according to our assumptions 
	and $\zeta'\equiv 0$ on $\R\backslash[-1,1]$, and we recall the definition 
	(\ref{eta}) of $\eta$ to find $c_6>0$
  satisfying
  \bas
	\io (x+\eps)^{\alpha+\beta-\gamma} \cdot (x+\eps)^{-2\beta} \psi_x^2(x) dx
	&=& \frac{1}{\eta^2} \io (x+\eps)^{\alpha-\beta-\gamma} \cdot \zeta \Big( \frac{x-x_0}{\eta}\Big) dx \\
	&\le& \frac{1}{\eta^2} \cdot (L+1)^{\alpha-\beta-\gamma} \cdot 2\eta \cdot \|\zeta'\|_{L^\infty(\R)} \\
	&=& c_6 M^{-\frac{1}{\theta}} (t_2-t_1)^{-\frac{1}{2\theta+3}} \\
	&\le& c_6 (t_2-t_1)^{-\frac{1}{2\theta+3}},
  \eas
  because of $M>1$ and thus $M^{-\frac{1}{\theta}}<1$. 
	Similarly, with some $c_7>0$ we have
  \bas
	\io (x+\eps)^{\alpha+\beta-\gamma} \cdot (x+\eps)^{-2\beta-2} \psi^2(x) dx
	&=& \io (x+\eps)^{\alpha-\beta-\gamma-2} \psi^2(x)dx \\
	&\le& (L+1)^{\alpha-\beta-\gamma-2} \cdot 2\eta \\
	&=& c_7 \cdot M^\frac{1}{\theta} (t_2-t_1)^\frac{1}{2\theta+3},
  \eas
  whence altogether
  \bas
	\bigg( \io (x+\eps)^{\alpha+\beta-\gamma} \Big[(x+\eps)^{-\beta} \psi(x) \Big]_x^2 dx \bigg)^\frac{1}{2}
	\le c_8 \Big\{ (t_2-t_1)^{-\frac{1}{2(2\theta+3)}} 
		+ M^\frac{1}{2\theta} (t_2-t_1)^\frac{1}{2(2\theta+3)} \Big\}
  \eas
  holds with some $c_8>0$. 
  Therefore, (\ref{14.5}), (\ref{14.6}), and (\ref{14.7}) in the limit $\delta\searrow 0$ yield $c_9>0$ such that
  \bas
	\frac{c_3}{4} M^{1+\frac{1}{\theta}} (t_2-t_1)^\frac{\theta+1}{2\theta+3}
	&\le& c_9 \Big\{ (t_2-t_1)^{\frac{1}{2}-\frac{1}{2(2\theta+3)}} 
		+ M^\frac{1}{2\theta} (t_2-t_1)^{\frac{1}{2}+\frac{1}{2(2\theta+3)}} \Big\} \\
	&=& c_9 \Big\{ (t_2-t_1)^\frac{\theta+1}{2\theta+3} 
		+ M^\frac{1}{2\theta} (t_2-t_1)^\frac{\theta+2}{2\theta+3} \Big\},
  \eas
  which implies the inequality
  \bas
	c_4 M^{1+\frac{1}{\theta}}
	&\le& c_9 \Big\{ 1 + M^\frac{1}{2\theta} (t_2-t_1)^\frac{1}{2\theta+3} \Big\} \\
	&\le& c_9 \Big\{ 1 + M^\frac{1}{2\theta} (T_0(A,B))^\frac{1}{2\theta+3} \Big\}.
  \eas
  Since $1+\frac{1}{\theta}>\frac{1}{2\theta}$, 
  this gives an upper bound for $M$ and thus yields the desired contradiction if $M$ has been chosen suitably large initially.
  This proves the H\"older estimate (\ref{14.3}) in time, and combining the latter with (\ref{14.2}) completes the proof.
\qed
\mysection{Proof of Theorem \ref{theo_final}}\label{sec.proof}
Let $u_0\in W^{1,2}_\gamma(\Omega)$, where $\gamma<1$, and let 
$(\eps_j)_{j\in\N} \subset (0,1)$ be a sequence satisfying $\eps_j\to 0$ as $j\to\infty$.
By a standard approximation
argument, we may construct a sequence of functions $(u_{0\eps_j})_{j\in\N}\subset 
C^\infty(\overline\Omega)$ such that 
\be{i1}
	u_{0\eps}>0 \quad \mbox{in } \bar\Omega
	\qquad \mbox{and} \qquad
	u_{0\eps x} \in C_0^\infty(\Omega)
	\qquad \qquad \mbox{for all } \eps \in (\eps_j)_{j\in\N},
\ee
and such that
\be{i2}
	u_{0\eps} \to u_0	
	\quad \mbox{in } W^{1,2}_\gamma(\Omega)
	\qquad \mbox{as } \eps=\eps_j\searrow 0.
\ee
The following lemma asserts that under the assumptions on $n$, $\alpha$, $\beta$, 
and $\gamma$ required in Lemma \ref{lem13},
the corresponding solutions of (\ref{0eps}) emanating from $u_{0\eps_j}$
have their maximal existence time bounded from below for all sufficiently small $\eps\in (\eps_j)_{j\in\N}$,
and moreover they accumulate at some continuous weak solution of (\ref{0}).
\begin{lem}\label{lem30}
  Let $n$, $\alpha$, $\beta$, and $\gamma$ be as in Theorem \ref{theo_final}.
  Then for all $A>0$ and $B>0$, there exists $T(A,B)\in (0,1)$ such that, whenever
  $u_0\in W^{1,2}_{loc}(\Omega)$ is nonnegative and satisfies
  \be{30.2_3}
	\io x^\gamma u_{0x}^2(x)dx \le A
  	\qquad \mbox{and} \qquad
	\io x^\beta u_0(x) dx \le B,
  \ee
  the following holds:
  For any $(\eps_j)_{j\in\N} \subset (0,\eps_0)$ such that $\eps_j\to 0$ as $j\to\infty$
  and each $(u_{0\eps_j})_{j\in\N} \subset C^\infty(\bar\Omega)$ 
  fulfilling (\ref{i1}) and (\ref{i2}),
  the problem (\ref{0eps}) possesses a unique positive classical solution $u_\eps\in C^{4,1}(\bar\Omega \times [0,T(A,B)])$
  for all sufficiently small $\eps\in (\eps_j)_{j\in\N}$, and moreover, there exists
  a subsequence $(\eps_{j_l})_{l\in\N}$ such that
  \be{30.7}
	u_\eps \to u
	\quad \mbox{in } C^0(\bar\Omega \times [0,T(A,B)])
	\qquad \mbox{as } \eps=\eps_{j_l}\to 0
  \ee
  with some continuous weak solution $u$ of (\ref{0}) in $\Omega \times (0,T(A,B))$.
\end{lem}
\proof
  We claim that the statement is valid if we let $T\equiv T(A,B):=T_0(A+1,B+1)$ 
  with $T_0$ as provided by Lemma \ref{lem13}.
  To verify this, we first note that, according to (\ref{30.2_3})
	and upon passing to subsequences,
  we may assume that
  \be{30.9}
	\io (x+\eps_j)^\gamma u_{0\eps_j x}^2(x) \le A+1
	\qquad \mbox{for all } j\in\N
  \ee
  and
  \be{30.10}
	\io (x+\eps_j)^\beta u_{0\eps_j}(x)dx \le B+1
	\qquad \mbox{for all } j\in\N.
  \ee
  For $\eps\in (\eps_j)_{j\in\N}$, we then let $u_\eps$ denote the corresponding positive classical solution of (\ref{0eps}),
  that is, of the initial-boundary value problem
  \bas
	\left\{ \begin{array}{l}
	u_{\eps t} = \frac{1}{(x+\eps)^\beta} \cdot \Big\{ -g_\eps(x) u_\eps^n u_{\eps xx}
		+2g_\eps(x) u_\eps^{n-1} u_{\eps x}^2 \Big\}_{xx}, \qquad x\in \Omega, \ t>0, \\[1mm]
	u_{\eps x}= u_{\eps xxx}=0, \qquad x\in\pO, \ t>0, \\[1mm]
	u_\eps(x,0)=u_{0\eps}(x), \qquad x\in\Omega,
	\end{array} \right.
  \eas
  which according to Lemma \ref{lem9} exists up to a maximal time $T_\eps \in (0,\infty]$ having the property stated
  in (\ref{9.1}).
  We divide the proof into four steps.\abs
  \underline{Step 1.} \quad
  We first show that actually $T_\eps \ge T$.\\
  To see this, we apply Lemma \ref{lem13} to find, upon passing to a
  subsequence if necessary, that for some $c_1>0$ and all 
	$\eps \in (\eps_j)_{j\in\N}$, we have
  \bea{30.11}
	&&\hspace*{-10mm}\sup_{t\in (0,\hat T_\eps)} \io (x+\eps)^\gamma u_{\eps x}^2(x,t)dx
	\nn \\
	&+& \int_0^{\hat T_\eps} \io (x+\eps)^{\alpha-\beta+\gamma} u_\eps^{n-2} u_{\eps x}^2 u_{\eps xx}^2dxdt
	+ \int_0^{\hat T_\eps} \io (x+\eps)^{\alpha-\beta+\gamma} u_\eps^{n-4} u_{\eps x}^6dxdt \nn\\
	&+& \int_0^{\hat T_\eps} \io (x+\eps)^{\alpha-\beta+\gamma-2} u_\eps^n u_{\eps xx}^2dxdt
	+ \int_0^{\hat T_\eps} \io (x+\eps)^{\alpha-\beta+\gamma-2} u_\eps^{n-2} u_{\eps x}^4dxdt \nn\\[2mm]
	&\le& c_1,
  \eea
  where $\hat T_\eps :=\min\{T_\eps,T\}$.
  Since
  \bas
	\sup_{t\in (0,\hat T_\eps)} \io (x+\eps)^\beta u_\eps(x,t)dx
	= \io (x+\eps)^\beta u_{0\eps}(x)dx  \le B+1
  \eas
  by Lemma \ref{lem8}, Lemma \ref{lem32} yields $c_2>0$ such that for all $\eps \in (\eps_j)_{j\in\N}$,
  \be{30.12}
	\sup_{t\in (0,\hat T_\eps)} \|u_\eps(\cdot,t)\|_{L^\infty(\Omega)} \le c_2.
  \ee
  Moreover, for fixed $\eps \in (\eps_j)_{j\in\N}$ we may apply Lemma \ref{lem31} with $\gamma$ replaced by $0$
  to see that (\ref{30.11}) entails that with some $c_3(\eps)>0$, the spatial H\"older estimate
  \be{30.13}
	|u_\eps(x,t)-u_\eps(y,t)| \le c_3(\eps) |x-y|^\frac{1}{2}
  \ee
  is valid for all $x,y\in\bar\Omega$ and any $t\in (0,\hat T_\eps)$. 
  Now assuming that $T_\eps<T$ for some $\eps\in (\eps_j)_{j\in\N}$, in view of the extensibility criterion in
  Lemma \ref{lem9} and the inequality (\ref{30.12}), we would have
  \bas
	u_\eps(x_k,t_k) \to 0
	\qquad \mbox{as } k\to\infty
  \eas
  with some $(x_k)_{k\in\N}\subset\Omega$ and $(t_k)_{k\in\N} \subset (0,T_\eps)$, where we may assume that
  $x_k\to x_0$ and $t_k \nearrow T_\eps$ as $k\to\infty$ with some $x_0\in\bar\Omega$.
  According to (\ref{30.12}), (\ref{30.13}), and the Arzel\`a-Ascoli theorem, we may pass to subsequences to
  achieve that with some $v\in C^0(\bar\Omega)$ we have
  \be{30.131}
	u_\eps(\cdot, t_k) \to v
	\quad \mbox{in } C^0(\bar\Omega)
	\qquad \mbox{as } k\to\infty,
  \ee
  and conclude that $v(x_0)=0$ and hence, again by (\ref{30.13}), that
  \be{30.132}
	0 \le v(x) \le c_3(\eps) |x-x_0|^\frac{1}{2}
	\qquad \mbox{for all } x\in\Omega.
  \ee
  This, however, contradicts the outcome of Lemma \ref{lem12}: The latter, namely, along with (\ref{30.12}) 
  implies that with some $c_4(\eps)>0$ we have
  \bas
	\io \frac{1}{u_\eps^2(x,t_k)} dx \le c_4(\eps)
	\qquad \mbox{for all } k\in\N,
  \eas
  so that Fatou's lemma and (\ref{30.132}) give
  \bas
	\frac{1}{c_3^2(\eps)}\io\frac{1}{|x-x_0|}dx\le
	\io \frac{1}{v^2(x)} dx \le c_4(\eps),
  \eas
  which is impossible.\abs
  \underline{Step 2.} \quad
  We next construct the limit function $u$.\\\
  To achieve this, we observe that since $T_\eps \ge T$ according to the above
	arguments, 
  we may replace $\hat T_\eps$ by $T$ in (\ref{30.11}) and (\ref{30.12})
  and apply Lemma \ref{lem14} to derive the $\eps$-independent estimate
  \bas
	\|u_\eps\|_{C^{\theta,\frac{\theta}{2\theta+3}}(\bar\Omega \times [0,T])} \le c_5
	\qquad \mbox{for all } \eps\in (\eps_j)_{j\in\N}
  \eas
  with a certain $c_5>0$. Therefore, the Arzel\`a-Ascoli theorem yields a subsequence, again denoted by 
  $(\eps_j)_{j\in\N}$, and a nonnegative function $u \in C^0(\bar\Omega \times [0,T])$ such that
  \be{30.18}
	u_\eps \to u
	\qquad \mbox{in } C^0(\bar\Omega \times [0,T])
  \ee
  as $\eps=\eps_j\searrow 0$. Moreover, interior parabolic regularity theory \cite{Fri64} shows that
  $(u_{\eps_j})_{j\in\N}$ is relatively compact in
  $C^{4,1}_{loc}(((0,L]\times (0,T]) \cap \{u>0\})$, and hence we may assume that as $\eps=\eps_j\searrow 0$,
  we also have
  \be{30.19}
	u_\eps \to u
	\qquad \mbox{in } C^{4,1}_{loc}(\pos), \qquad
	\mbox{where } \pos:=\Big((0,L]\times (0,T]\Big) \cap \{u>0\}.
  \ee
  \underline{Step 3.} \quad
  We proceed to verify that there exists a null set $N\subset (0,T)$ such that for all $t\in (0,T)\setminus N$,
  $u(\cdot,t)$ is differentiable at $x=L$ with $u_x(L,t)=0$ for all $t\in (0,T)\setminus N$.\\
  To this end, we note that (\ref{30.11}) in particular implies that for some $c_6>0$, we have
  \bas
	\int_0^T \int_\frac{L}{2}^L u_\eps^{n-2} u_{\eps x}^2 u_{\eps xx}^2dxdt
	+ \int_0^T \int_\frac{L}{2}^L u_\eps^{n-4} u_{\eps x}^6 dxdt
	\le c_6
	\qquad \mbox{for all } \eps\in (\eps_j)_{j\in\N}
  \eas
  and since 
  \bas
	\bigg\{ \Big(u_\eps^\frac{n+2}{4}\Big)_x^2 \bigg\}_x^2
	&=& \Big(\frac{n+2}{4}\Big)^4 \cdot \Big\{ 2u_\eps^\frac{n-2}{2} u_{\eps x} u_{\eps xx}
	+ \frac{n-2}{2} u_\eps^\frac{n-4}{2} u_{\eps x}^3 \Big\}^2 \\
	&\le& 2\Big(\frac{n+2}{4}\Big)^4 \cdot \Big\{ 4u_\eps^{n-2} u_{\eps x}^2 u_{\eps xx}^2
	+ \Big(\frac{n-2}{2}\Big)^2 u_\eps^{n-4} u_{\eps x}^6 \Big\}
	\qquad \mbox{in } \Omega \times (0,T),
  \eas
  we thus find $c_7>0$ fulfilling
  \be{30.89}
	\int_0^T \int_\frac{L}{2}^L \bigg\{ \Big(u_\eps^\frac{n+2}{4}\Big)_x^2\bigg\}_x^2dxdt \le c_7
	\qquad \mbox{for all } \eps\in (\eps_j)_{j\in\N}.
  \ee
  Since $u_{\eps x}(L,t)=0$ and hence $(u_\eps^\frac{n+2}{4})_x(L,t)=0$ for all $t\in (0,T)$ by (\ref{0eps}),
  using the Cauchy-Schwarz inequality, we obtain
  \bea{30.90}
	\Big(u_\eps^\frac{n+2}{4}\Big)_x^2(x,t) 
	&=& - \int_x^L \bigg\{ \Big(u_\eps^\frac{n+2}{4}\Big)_x^2 \bigg\}_x (y,t) dy \nn\\
	&\le& (L-x)^\frac{1}{2} \cdot a_\eps(t)
	\qquad \mbox{for all } x\in \Big(\frac{L}{2},L\Big) \mbox{ and } t\in (0,T),
  \eea
  where
  \bas
	a_\eps(t):=\int_\frac{L}{2}^L \bigg\{ \Big(u_\eps^\frac{n+2}{4}\Big)_x^2\bigg\}_x^2(y,t)dy,
	\qquad t\in (0,T).
  \eas
  Again by the Cauchy-Schwarz inequality, (\ref{30.90}) in turn implies that
  \bas
	\Big| u_\eps^\frac{n+2}{4}(L-x,-t)-u_\eps^\frac{n+2}{4}(x,t)\Big|
	&=& \bigg| \int_x^L \Big(u_\eps^\frac{n+2}{4}\Big)_x(y,t) dy \bigg| \\
	&\le& (L-x)^\frac{1}{2} \cdot \bigg\{ \int_x^L (L-y)^\frac{1}{2} \cdot a_\eps^\frac{1}{2}(t) dy 
		\bigg\}^\frac{1}{2} \\
	&=& \sqrt{\frac{2}{3}} (L-x)^\frac{5}{4} \cdot a_\eps^\frac{1}{4}(t)
	\qquad \mbox{for all } x\in \Big(\frac{L}{2},L\Big) \mbox{ and } t\in (0,T),
  \eas
  by (\ref{30.89}) and the definition of $a_\eps$ meaning that
  \bas
	\int_0^T \sup_{x\in (\frac{L}{2},L)} 
	\frac{\Big|u_\eps^\frac{n+2}{4}(L,t)-u_\eps^\frac{n+2}{4}(x,t)\Big|^4}{(L-x)^5} \, dt
	\le \frac{4}{9} c_7.
  \eas
  Using (\ref{30.18}) and Fatou's lemma, from this we conclude that
  \bas
	\int_0^T \sup_{x\in (\frac{L}{2},L)} 
	\frac{\Big|u^\frac{n+2}{4}(L,t)-u^\frac{n+2}{4}(x,t)\Big|^4}{(L-x)^5} \, dt
	\le \frac{4}{9} c_7,
  \eas
  so that in particular we can find a null set $N\subset (0,T)$ such that for all $t\in (0,T)\setminus N$,
  \bas
	b(t):=\sup_{x\in (\frac{L}{2},L)}  
		\frac{\Big|u^\frac{n+2}{4}(L,t)-u^\frac{n+2}{4}(x,t)\Big|^4}{(L-x)^5}
  \eas
  is finite.\\
  Now if $t\in (0,T)\setminus N$ is such that $u(L,t)>0$, then from (\ref{30.19}) we clearly infer the existence
  of $u_x(L,t)=\lim_{\eps=\eps_j\searrow 0} u_{\eps x}(L,t)=0$.
  On the other hand, if $t\in (0,T)\setminus N$ is such that $u(L,t)=0$, then according to the definition of $b(t)$,
  we obtain
  \bas
	u^{n+2}(x,t) \le b(t) \cdot (L-x)^5
	\qquad \mbox{for all } x\in \Big(\frac{L}{2},L\Big),
  \eas
  that is,
  \bas
	\bigg| \frac{u(L,t)-u(x,t)}{L-x} \bigg|
	&=& \frac{u(x,t)}{L-x} \\
	&\le& \frac{\Big\{ b(t)\cdot (L-x)^5\Big\}^\frac{1}{n+2}}{L-x} \\
	&=& b^\frac{1}{n+2}(t) \cdot (L-x)^\frac{3-n}{n+2}
	\qquad \mbox{for all }  x\in \Big(\frac{L}{2},L\Big),
  \eas
  so that, since $n<3$, we infer that also in this case $u_x(L,t)$ exists and vanishes.\abs
  \underline{Step 4. } \quad
  We finally show that $u$ furthermore satisfies the integral identity (\ref{0w}).\\
  To prepare this, let us first derive two further estimates from (\ref{30.11}):
  Namely, for any $q\in [1,2)$, (\ref{30.11}) along with Lemma \ref{lem1} and the H\"older inequality
  implies that for all measurable $\Omega_0\subset\Omega$ and any measurable $Q \subset \Omega_0 \times (0,T)$, 
  we have
  \bea{30.20}
	\iint_Q \Big|g_\eps(x) u_\eps^n u_{\eps xx} \Big|^q dxdt
	&\le& \!\!\iint_Q (x+\eps)^{q\alpha} u_\eps^{nq} |u_{\eps xx}|^qdxdt \nn\\
	&\le& \!\!\Big(\iint_Q (x+\eps)^{\alpha-\beta+\gamma-2} u_\eps^n u_{\eps xx}^2 dxdt\Big)^\frac{q}{2}
	\Big( \iint_Q (x+\eps)^\frac{q(\alpha+\beta-\gamma+2)}{2-q} u_\eps^\frac{nq}{2-q}dxdt \Big)^\frac{2-q}{2}
		\nn\\
	&\le& \!\! c_1^\frac{q}{2} T^\frac{2-q}{2} \cdot 
		\Big( \int_{\Omega_0} (x+\eps)^\frac{q(\alpha+\beta-\gamma+2)}{2-q} dx\Big)^\frac{2-q}{2}
		\cdot \|u_\eps\|_{L^\infty(Q)}^\frac{nq}{2},
  \eea
  whereas similarly,
  \bea{30.21}
	\iint_Q \Big|g_\eps(x) u_\eps^{n-1} u_{\eps x}^2 \Big|^qdxdt
	&\le& \iint_Q (x+\eps)^{q\alpha} u_\eps^{(n-1)q} |u_{\eps x}|^{2q}dxdt \nn\\
	&\le& \Big( \iint_Q (x+\eps)^{\alpha-\beta+\gamma-2} u_\eps^{n-2} u_{\eps x}^4 dxdt\Big)^\frac{q}{2} \nn \\
	&&{}\times \Big( \iint_Q (x+\eps)^\frac{q(\alpha+\beta-\gamma+2)}{2-q} u_\eps^\frac{nq}{2-q}dxdt \Big)^\frac{2-q}{2} \nn\\
	&\le& c_1^\frac{q}{2} T^\frac{2-q}{2} \cdot
	\Big( \int_{\Omega_0} (x+\eps)^\frac{q(\alpha+\beta-\gamma+2)}{2-q}dx \Big)^\frac{2-q}{2}
	\cdot \|u_\eps\|_{L^\infty(Q)}^\frac{nq}{2}.
  \eea
	In view of our assumptions $\alpha>3$, $\beta>-1$, and $\gamma<1$,
	we infer that $\alpha+\beta-\gamma+2>3$.
  Hence, picking any $q\in (1,2)$ we know that
  \bas
	c_8:=\sup_{\eps \in (0,1)} \Big( \io (x+\eps)^\frac{q(\alpha+\beta-\gamma+2)}{2-q}dx \Big)^\frac{2-q}{2}
  \eas
  is finite. Then (\ref{30.20}) and (\ref{30.21}), applied to $Q:=Q_T:=\Omega\times (0,T)$, 
  show that because of $q>1$, we may
  pass to a further subsequence to achieve that
  \bas
	g_\eps(x) u_\eps^n u_{\eps xx} \wto w
	\qquad \mbox{in } L^q(\Omega\times (0,T))
  \eas
  and
  \bas
	g_\eps(x) u_\eps^{n-1} u_{\eps x}^2 \wto z
	\qquad \mbox{in } L^q(\Omega\times (0,T))
  \eas
  as $\eps=\eps_j\to 0$ with some $w$ and $z$ belonging to $L^q(\Omega\times (0,T))$. 
  In view of the pointwise convergence properties $u_{\eps x} \to u_x$ and $u_{\eps xx} \to u_{xx}$
  inside $\pos$, as guaranteed by (\ref{30.19}), we may identify these limits to obtain that actually
  \be{30.22}
	g_\eps(x) u_\eps^n u_{\eps xx} \wto x^\alpha u^n u_{xx}
	\qquad \mbox{in } L^q(\pos)
  \ee
  and
  \be{30.23}
	g_\eps(x) u_\eps^{n-1} u_{\eps x}^2 \wto x^\alpha u^{n-1} u_x^2
	\qquad \mbox{in } L^q(\pos),
  \ee
  because $n>1$. \\
  Next, outside the set $\pos$, we may use that $u_\eps \to 0$ uniformly in $Q_T \setminus \pos$ to infer upon 
  another application of (\ref{30.19}) and (\ref{30.20}) to $q:=1$ that
  \be{30.24}
	\iint_{Q_T \setminus \pos} g_\eps(x) u_\eps^n |u_{\eps xx}|dxdt \to 0
  \ee
  and 
  \be{30.25}
	\iint_{Q_T \setminus \pos} g_\eps(x) u_\eps^{n-1} u_{\eps x}^2dxdt \to 0
  \ee
  as $\eps=\eps_j\to 0$, noting here again that $\alpha+\beta-\gamma+2>0$ by assumption.\abs
  Now for the verification of (\ref{0w}), we fix any $\phi \in C_0^\infty(\bar\Omega \times [0,T))$
  such that $\phi_x=0$ at $x=L$. 
  We then approximate $\phi$ by letting
  \be{30.255}
	\phi_\delta(x,t) := \phi(0,t) + \int_0^x \zeta_\delta(y) \phi_x(y,t) dy,
	\qquad (x,t) \in\bar\Omega \times [0,T),
  \ee 
  for $\delta \in (0,\frac{L}{2})$, where
  \bas
	\zeta_\delta(x):=\zeta \Big(\frac{x}{\delta}\Big), \qquad x\in\bar\Omega,
  \eas
  with a fixed cut-off function $\zeta \in C^\infty\R)$ such that
  $\zeta\equiv 0$ in $(-\infty,1]$, $\zeta\equiv 1$ in $[2,\infty)$ and $0 \le \zeta' \le 2$ on $\R$.
  This construction ensures that $\phi_{\delta x}$ vanishes at both $x=L$ and $x=0$, so that upon
  multiplying (\ref{0eps}) by $(x+\eps)^\beta \phi_{\delta}$, we may integrate by parts, again using
  Lemma \ref{lem_J}, to obtain
  \bea{30.26}
	- \int_0^T \io (x+\eps)^\beta u_\eps \phi_{\delta t}dxdt
	&-& \io (x+\eps)^\beta u_{0\eps}(x) \phi_\delta(x,0)dx \nn\\
	&=& \int_0^T \io \Big[ - g_\eps(x) u_\eps^n u_{\eps xx} + 2g_\eps(x) u_\eps^{n-1} u_{\eps x}^2 \Big]
		\cdot \phi_{\delta xx}dxdt
  \eea
  for all $\eps \in (\eps_j)_{j\in\N}$ and all $\delta \in (0,\frac{L}{2})$. Here, from (\ref{30.18}) and
  the fact that $u_{0\eps} \to u_0$ in $C^0(\bar\Omega)$ by (\ref{i2}) and the restriction $\gamma<1$, it is clear that
  \bas
	- \int_0^T \io (x+\eps)^\beta u_\eps \phi_{\delta t}dxdt \to - \int_0^T \io x^\beta u \phi_t dxdt
  \eas
  and
  \bas
	- \io (x+\eps)^\beta u_{0\eps}(x) \phi_\delta(x,0)dx
	\to - \io x^\beta u_0(x) \phi_\delta(x,0)dx
  \eas
  as $\eps=\eps_j\to 0$, whereas (\ref{30.22})-(\ref{30.25}) warrant that
  \bas
	\int_0^T \io \Big[ - g_\eps(x) u_\eps^n u_{\eps xx} + 2g_\eps(x) u_\eps^{n-1} u_{\eps x}^2 \Big]
		\cdot \phi_{\delta xx}dxdt
	\to \iint_{\pos} [-x^\alpha u^n u_{xx} + 2x^\alpha u^{n-1} u_x^2 ] \cdot \phi_{\delta xx}dxdt
  \eas
  as $\eps=\eps_j\to 0$; hence, (\ref{30.26}) yields
  \be{30.27}
	- \int_0^T \io x^\beta u\phi_{\delta t}dxdt
	- \io x^\beta u_0(x)\phi_\delta(x,0)dx
	= \iint_{\pos} [-x^\alpha u^n u_{xx} + 2x^\alpha u^{n-1} u_x^2 ] \cdot \phi_{\delta xx}dxdt
  \ee
  for all $\delta \in (0,\frac{L}{2})$. Now, taking $\delta\searrow 0$, we observe that by (\ref{30.255}),
  \bas
	\phi_{\delta xx}(x,t) = \zeta_\delta(x) \cdot \phi_{xx}(x,t)
	+ \frac{1}{\delta} \zeta' \Big(\frac{x}{\delta}\Big) \cdot \phi_x(x,t)
	\qquad \mbox{for all } (x,t)\in \Omega\times (0,T),
  \eas
  so that since $0\le \zeta' \le 2$ we find that
  \bas
	& & \hspace*{-10mm}
	\bigg| \iint_{\pos} [-x^\alpha u^n u_{xx} + 2x^\alpha u^{n-1} u_x^2 ] \cdot \phi_{\delta xx}dxdt
	- \iint_{\pos} [-x^\alpha u^n u_{xx} + 2x^\alpha u^{n-1} u_x^2 ] \cdot \phi_{xx}dxdt \bigg|  \\
	&\le& \|\phi_{xx}\|_{L^\infty(\Omega \times (0,T))} \cdot
	\bigg| \iint_{\pos} [-x^\alpha u^n u_{xx} + 2x^\alpha u^{n-1} u_x^2 ] \cdot (1-\zeta_\delta(x))dxdt \bigg| \\
	& & + \frac{2}{\delta} \|\phi_x\|_{L^\infty(\Omega\times(0,T))} \cdot \iint_{S_\delta} x^\alpha u^n |u_{xx}|dxdt
	+ \frac{4}{\delta} \|\phi_x\|_{L^\infty(\Omega\times(0,T))} \cdot \iint_{S_\delta} x^\alpha u^{n-1} u_x^2 dxdt \\[2mm]
	&=:& I_1(\delta)+I_2(\delta)+I_3(\delta),
  \eas
  where $S_\delta:=((0,2\delta)\times (0,T)) \cap \pos$.
  Clearly, 
  \bas
	I_1(\delta) \to 0
	\qquad \mbox{as } \delta\searrow 0
  \eas
  by the dominated convergence theorem in conjunction with the integrability property
  $-x^\alpha u^n u_{xx}+2x^\alpha u^{n-1} u_x^2 \in L^1(\pos)$ asserted by (\ref{30.22}) and (\ref{30.23}).
  Moreover, applying (\ref{30.19}) and (\ref{30.20}) to $Q:=((0,2\delta)\times (0,T)) \cap \pos$
  and $q:=1$ and once more recalling (\ref{30.22}) and (\ref{30.23}), we see that
  \bas
	\iint_{S_\delta} x^\alpha u^n |u_{xx}|dxdt
	\le c_9 \Big( \int_0^{2\delta} x^{\alpha+\beta-\gamma+2} dx \Big)^\frac{1}{2} 
	\le c_{10} \delta^\frac{\alpha+\beta-\gamma+3}{2}
  \eas
  and similarly 
  \bas
	\iint_{S_\delta} x^\alpha u^{n-1} u_x^2dxdt
	\le c_{11} \delta^\frac{\alpha+\beta-\gamma+3}{2}
  \eas
  for all $\delta \in (0,\frac{L}{2})$ with positive constants $c_9$, $c_{10}$, and $c_{11}$.
  As our hypotheses $\alpha>3$, $\beta>-1$, and $\gamma<1$ guarantee that 
  $\frac{\alpha+\beta-\gamma+3}{2}>1$, we thus obtain that also
  \bas
	I_2(\delta) + I_3(\delta) \to 0
	\qquad \mbox{as } \delta\searrow 0,
  \eas
  so that, since clearly $\phi_\delta \to \phi$ and $\phi_{\delta t} \to \phi_t$ uniformly in
  $\Omega \times (0,t)$, we conclude from (\ref{30.27})
	that indeed (\ref{0w}) is valid.
\qed
We can now prove our main result.\abs	
\proofc of Theorem \ref{theo_final}. \quad
  According to Lemma \ref{lem30} with $K>0$ as in (\ref{7.1}), we know
  that there exists $T>0$ and a continuous weak solution $u$ of (\ref{0}) in $\Omega \times (0,T)$, 
  which due to Lemma \ref{lem7} and the approximation statement in Lemma \ref{lem30} has
  the additional regularity property
  \be{34.33}
	u\in L^\infty ((0,T);W^{1,2}_\gamma(\Omega))
  \ee
  and satisfies
  \be{34.3}
	\io x^\gamma u_x^2(x,t)dx \le \io x^\gamma u_{0x}^2(x)dx + K \int_0^t \io x^{\alpha-\beta+\gamma-6} u^{n+2}dxds
	\qquad \mbox{for a.e.~} t\in (0,T).
  \ee
  From Lemma \ref{lem8} combined with Lemma \ref{lem30}, we infer that moreover
  \be{34.4}
	\io x^\beta u(x,t)dx = B_0:= \io x^\beta u_0(x)dx
	\qquad \mbox{for all } t\in (0,T).
  \ee
  Therefore, 
  \bas
	\tm:= \sup \Big\{ T>0 \ &\Big|& \ 
	\mbox{There exists a continuous weak solution $u$ of (\ref{0}) in $\Omega\times (0,T)$} \\
	& & \mbox{which satisfies (\ref{34.33}), (\ref{34.3}) and (\ref{34.4})} \ \Big\}
	\quad \le\infty
  \eas
  is well-defined, and it remains to show that (\ref{extend}) holds.\\
  Indeed, let us assume on the contrary that $\tm<\infty$ but $u\le M$ in $\Omega\times (0,\tm)$ for some $M>0$.
  Then (\ref{34.3}) would imply that
  \bas
	\io x^\gamma u_x^2(x,t)dx
	\le A_0 :=\io x^\gamma u_{0x}^2(x)dx
	+ K M^{n+2} \tm \io x^{\alpha-\beta+\gamma-6}dx
	\quad \mbox{for a.e.~} t\in (0,\tm),
  \eas
  where our assumption $\gamma>5-\alpha+\beta$ ensures that $\alpha-\beta+\gamma-6>-1$ and hence $A_0<\infty$.
  We could thus pick some $t_0\in (0,\tm)$ such that
  \bas
	t_0>\tm-\frac{1}{2} T(A_0,B_0)
	\qquad \mbox{and} \qquad
	\io x^\gamma u_x^2(x,t_0)dx \le A_0,
  \eas
  to see upon another application of Lemma \ref{lem30} to $A:=A_0$, $B:=B_0$ and
  \bas
	v_0(x):=u(x,t_0), \qquad x\in\Omega,
  \eas
  that the problem
  \bas
	\left\{ \begin{array}{l}
	v_t=\frac{1}{x^\beta} \cdot \Big\{ x^\alpha [-v^n v_{xx} + 2v^{n-1} v_x^2] \Big\}_{xx},
	\qquad x\in \Omega, \ t>0, \\[1mm]
	x^\alpha [-v^n u_{xx} + 2 v^{n-1} u_x^2]=x^\alpha [-v^n v_{xx} + 2 v^{n-1} v_x^2]_x=0, \qquad x=0, \ t>0, \\[1mm]
	v_x=v_{xxx}=0, \qquad x=L, \ t>0, \\[1mm]
	v(x,0)=v_0(x), \qquad x\in \Omega,
	\end{array} \right.
  \eas
  would possess a continuous weak solution $v$ in $\Omega \times (0,T(A_0,B_0))$ which, 
  again by Lemma \ref{lem7}, Lemma \ref{lem8}, and Lemma \ref{lem30}, would satisfy
  $v\in L^\infty((0,T(A_0,B_0));W^{1,2}_\gamma(\Omega))$ and
  \bas
	\io x^\gamma v_x^2(x,t) dx \le \io x^\gamma v_{0x}^2(x)dx + K \int_0^t \io x^{\alpha-\beta+\gamma-6} v^{n+2}dxds
	\quad \mbox{for a.e.~} t\in (0,T(A_0,B_0))	
  \eas
  as well as
  \bas
	\io x^\beta v(x,t) dx =B_0
	\qquad \mbox{for all } t\in (0,T(A_0,B_0)).
  \eas
  It can therefore easily be checked that
  \bas
	\tilde u(x,t):=\left\{ \begin{array}{ll}
	u(x,t) \qquad & \mbox{if } x\in \Omega \mbox{ and } t\in (0,t_0), \\[1mm]
	v(x,t-t_0) \qquad &\mbox{if } x\in\Omega \mbox{ and } t\in [t_0,t_0+T(A_0,B_0)),
	\end{array} \right.
  \eas
  would define a continuous weak solution $\tilde u$ of (\ref{0}) in 
	$\Omega\times (0,t_0+T(A_0,B_0))$, 
  yet fulfilling (\ref{34.33}), (\ref{34.3}), and (\ref{34.4}).
  As $t_0+T(A_0,B_0)>\tm$, this contradicts the definition of $\tm$.
\qed
{\bf Acknowledgment.} \quad
The first author acknowledges partial support from   
the Austrian Science Fund (FWF), grants P20214, P22108, I395, and W1245.
He thanks Miquel Escobedo for the hint on the paper \cite{JPR06},
which initiated this study.
\end{document}